\documentstyle[11pt]{article}

\newcommand{\barl}{ \overline }

\newcommand{\A}{ {\cal A} }
\newcommand{\B}{ {\cal B} }
\newcommand{\C}{ \mbox{\bf C} }
\newcommand{\D}{ {\cal D} }
\newcommand{\M}{ {\cal M} }
\newcommand{\MM}{ M_{2} ( \C I ) }

\newcommand{\R}{ \mbox{\bf R} }

\newcommand{\X}{ {\cal X} }
\newcommand{\Y}{ {\cal Y} }
\newcommand{\Z}{ \mbox{\bf Z} }

\newcommand{\ff}{ \varphi }
\newcommand{\expb}{ E_{\M} }
\newcommand{\expd}{ E_{\D} }
\newcommand{\etaz}{ \eta_{0} }
\newcommand{\Xalt}{ \X_{alt}^{*} }

\newcommand{\ncps}{ ( {\cal A} , \varphi ) }

\newcommand{\ncpsd}{ ( M_{d} ( \A ) , \varphi_{d} ) }
\newcommand{\ncpsm}{ ( {\cal M} , \psi ) }
\newcommand{\ncpst}{ ( M_{2} ( \A ) , \varphi_{2} ) }

\newcommand{\laphi}{ L^{2}  \ncps  }

\newcommand{\ldaphi}{ L^{2}  \ncpsd  }
\newcommand{\laphitwo}{ L^{2}  \ncpst  }
\newcommand{\lphi}{ L^{2} ( \varphi ) }

\newcommand{\lphid}{ L^{2} ( \varphi_{d} ) }
\newcommand{\lphitwo}{ L^{2} ( \varphi_{2} ) }

\newcommand{\ais}{ \{ a_{i} \}_{i \in I} }

\newcommand{\xiis}{ \{ \xi_{i} \}_{i \in I} }
\newcommand{\etais}{ \{ \eta_{i} \}_{i \in I} }
\newcommand{\aijs}{ \{ a_{ij} \}_{1 \leq i,j \leq d} }
\newcommand{\bijs}{ \{ b_{ij} \}_{1 \leq i,j \leq d} }

\newcommand{\vijs}{ \{ v_{ij} \}_{1 \leq i,j \leq d} }

\newcommand{\aijaijs}{ \{ a_{ij} , a_{ij}^{*} \}_{1 \leq i,j \leq d} }
\newcommand{\cijcijs}{ \{ c_{ij} , c_{ij}^{*} \}_{1 \leq i,j \leq d} }
\newcommand{\xiijxiijs}{ \{ \xi_{ij}, \xi_{ij}^{*} \}_{1 \leq i,j \leq d} }

\newcommand{\mataij}{ ( a_{ij} )_{i,j=1}^{d} }
\newcommand{\matbij}{ ( b_{ij} )_{i,j=1}^{d} }
\newcommand{\matcij}{ ( c_{ij} )_{i,j=1}^{d} }
\newcommand{\Entry}{\mbox{Entry}}

\begin{document}

\title{\bf Some minimization problems for the free analogue of 
the Fisher information}

\author{
Alexandru Nica 
\thanks{Research supported by a grant from the Natural Sciences 
and Engineering Research Council, Canada. } \\
Department of Pure Mathematics \\
University of Waterloo \\
Waterloo, Ontario N2L 3G1, Canada \\
anica@math.uwaterloo.ca 
\and
Dimitri Shlyakhtenko
\thanks{Supported in part by a Graduate Research Fellowship of the
National Science Foundation, USA. } \\
Department of Mathematics \\
University of California at Los Angeles \\
Los Angeles, CA 90095-1555, U.S.A. \\
shlyakht@member.ams.org
\and
Roland Speicher
\thanks{Supported by a Heisenberg Fellowship of the DFG, Germany. } \\
Institut f\"ur Angewandte Mathematik \\
Universit\"at Heidelberg \\
D-69120 Heidelberg, Germany \\
roland.speicher@urz.uni-heidelberg.de } 

\date{ }

\maketitle

\begin{abstract}
We consider the free non-commutative analogue $\Phi^{*}$, 
introduced by D. Voiculescu, 
\setcounter{page}{0}
of the concept of Fisher information for random variables. We determine 
the minimal possible value of $\Phi^{*} ( a,a^{*} )$, if $a$ is a 
non-commutative random variable subject to the constraint that the 
distribution of $a^{*}a$ is prescribed. More generally, we obtain the 
minimal possible value of $\Phi^{*} ( \ \aijaijs \ )$, if $\aijs$ is a 
family of non-commutative random variables such that the distribution 
of $A^{*}A$ is prescribed, where $A$ is the matrix $\mataij$. The
$d \times d$-generalization is obtained from the case $d=1$ via a result 
of independent interest, concerning the minimal value of 
$\Phi^{*} ( \aijaijs )$ when the matrix $A = \mataij$ and its adjoint 
have a given joint distribution. (A version of this result describes 
the minimal value of $\Phi^{*} ( \bijs )$ when the matrix $B = \matbij$ 
is selfadjoint and has a given distribution.)

We then show how the minimization results obtained for $\Phi^{*}$
lead to maximization results concerning the free entropy $\chi^{*}$,
also defined by Voiculescu.
\end{abstract}

\newpage

\setcounter{section}{1}
\setcounter{equation}{0}
{\large\bf 1. Introduction} 

$\ $

In this paper we determine the minimal possible value of the free
Fisher information $\Phi^{*} ( a,a^{*} )$, if $a$ is a non-commutative 
random variable subject to the constraint that the distribution of 
$a^{*}a$ is prescribed. More generally, we obtain the minimal possible 
value of $\Phi^{*} ( \ \aijaijs \ )$, if $\aijs$ is a family of 
non-commutative random variables such that the distribution of 
$A^{*}A$ is prescribed, where $A$ is the matrix $\mataij$. The 
$d \times d$-generalization is obtained via a result of independent 
interest on the minimal free Fisher information of a family of matrix 
entries, when the distribution/$*$-distribution of the matrix itself is
given.

The framework we will consider is the one of a $W^{*}$-probability 
space $\ncps$, with $\varphi$ a faithful trace (i.e. -- $\A$ is a 
$W^{*}$-algebra, and $\varphi : \A \rightarrow \C$ is a normal faithful 
trace-state). An element $a \in \A$ will be referred to
as a ``non-commutative random variable'', and $\varphi (a)$ will be
called ``the expectation of $a$''. If $a = a^{*} \in \A$, then the 
unique probability measure with compact support $\mu$ on $\R$ which
has $\int_{- \infty}^{\infty} t^{n} \ d \mu (t)$ = $\varphi (a^{n})$, 
$\forall n \geq 0$, is called the distribution of $a$. An element 
$a = a^{*} \in \A$ is said to be semicircular of radius $r > 0$ if
its distribution is absolutely continuous with respect to the 
Lebesgue measure, with density $\rho (t) = 2( \pi r^{2} )^{-1}
\sqrt{ r^{2} - t^{2} }$ on $[ -r , r ]$.

A fundamental concept used throughout the paper is the one of freeness
for a family of subsets of $\A$. For the definition  and basic properties
of freeness, we refer the reader to \cite{VDN}, Chapter 2.

The free analogues of entropy and of Fisher information for random 
variables were introduced and studied in a series of papers of
D. Voiculescu (\cite{V-P1} -- \cite{V-P5}), in connection to the
isomorphism problem for the von Neumann algebras associated to free 
groups. Free analogues for some well-known inequalities concerning the 
Fisher information were obtained in this way. In particular, one has a
``free Cramer-Rao inequality'', which says the following: if 
$( x_{1} , \ldots , x_{n} )$ is an $n$-tuple of selfadjoint elements of
$\A$ such that the total variance $\varphi ( x_{1}^{2} + \cdots +
x_{n}^{2} )$ is prescribed, then the free Fisher information 
$\Phi^{*} ( x_{1} , \ldots , x_{n} )$ is minimized when the $x_{j}$'s
are semicircular of equal radii, and free (see \cite{V-P5}, Proposition
6.9). In the particular case $n=2$, if one sets $a := x_{1} +ix_{2}$ and 
works with $a,a^{*}$ instead of $x_{1},x_{2}$, then the free 
Cramer-Rao inequality can also be formulated like this: let $a$ be a 
non-commutative random variable, such that the expectation of $a^{*}a$
is prescribed; then the free Fisher information $\Phi^{*} (a,a^{*})$ is
minimized when $a$ is a circular element (which means, by definition,
that the real and imaginary part of $a$ are free and have semicircular 
distributions of equal radii).

In the present paper we examine a similar minimization problem, where not
only the expectation, but the whole distribution (i.e. the moments of all
orders) of $a^{*}a$ are prescribed. More precisely: given a probability
measure $\nu$ with compact support on $[0, \infty )$, 
what can be said about
\begin{equation}
\inf \{ \Phi^{*} ( a, a^{*} ) \ | \ a^{*}a \mbox{ has distribution }
\nu \} \ ?
\end{equation}
One cannot of course hope to have the infimum in (1.1) achieved by a
circular element; this is simply because, given $\nu$ as in (1.1), there
does not exist in general a circular element $a$ such that $a^{*}a$ has
distribution $\nu$. (In fact: if $a$ is circular, then the distribution 
of $a^{*}a$ can only be of the form 
$2( \alpha \pi )^{-1} \sqrt{ ( \alpha - t )/t } dt$ on $[ 0 , \alpha ]$ 
for some $\alpha > 0$ -- see \cite{VDN}, Section 5.1.)

A remarkable family of relatives of the circular element is provided by
the so-called ``$R$-diagonal elements'', introduced in \cite{NS}. There 
are several possible descriptions for the fact that an element
$a \in \A$ is $R$-diagonal. The one taken as starting point in
\cite{NS} is that the $R$-transform -- i.e. free analogue for the log of 
the Fourier transform -- of the pair $(a,a^{*})$ has a special form, which 
is in a certain sense ``diagonal''; this is in fact where the name of 
``$R$-diagonal'' comes from. In the present paper we will use an 
equivalent characterization of $R$-diagonality, described as follows:
$a$ is $R$-diagonal if and only if
the $*$-distribution of $a$ (i.e. the family of expectations of words
in $a$ and $a^{*}$) coincides with the $*$-distribution of an element of 
the form $up$, where $u$ is a unitary distributed according to
the Haar measure on the circle, $p=p^{*}$, and $\{ u,u^{*} \}$ is free 
from $\{ p \}$. The equivalence between the two characterizations of an
$R$-diagonal element is shown in \cite{NS}. The circular element is 
$R$-diagonal, e.g. because its polar decomposition is known to be of the
form $up$, with $u$ Haar unitary such that $\{ u,u^{*} \}$ is free from 
$\{ p \}$ (see \cite{VDN}, Section 5.1).

Now, given a probability measure $\nu$, with compact support on
$[ 0, \infty )$, there always exists an $R$-diagonal element $a$ such
that $a^{*}a$ has distribution $\nu$. This $a$ is ``unique up to 
isomorphism'', in the sense that the $*$-distribution of $a$ is 
completely determined (which in turn determines the
unital $W^{*}$-algebra generated by $a$); see Remark 3.3 below. The 
result we obtain is that the $R$-diagonal element attains the infimum
considered in (1.1). Moreover, finding the actual value of the infimum
is reduced to the calculation of a free Fisher information $\Phi^{*} 
( \mu )$, where $\mu$ is a symmetric distribution naturally associated 
to $\nu$; and for $\Phi^{*} ( \mu )$ one can use an explicit formula 
established in \cite{V-P1}. To summarize, we have:

$\ $

{\bf 1.1 Theorem.} Let $\nu$ be a probability measure with compact
support on $[ 0, \infty )$. Let $\mu$ be the symmetric probability measure
on $\R$ determined by the fact that $\mu (S) = \nu ( S^{2} )$ for every 
symmetric Borel set $S \subseteq \R$. Then
\begin{equation}
\min \{ \Phi^{*} ( a, a^{*} ) \ | \ a^{*}a \mbox{ has distribution }
\nu \} \ = \ 2 \Phi^{*} ( \mu ),
\end{equation}
and the minimum is attained when $a$ is $R$-diagonal. If in particular 
$\nu$ is absolutely continuous, with density $\rho$, then the quantities 
in (1.2) equal:
\begin{equation}
\frac{4}{3} \cdot \int_{0}^{\infty} t \rho (t)^{3} \ dt \ \in
\ [ 0, \infty ] .
\end{equation}

$\ $

The facts stated in Theorem 1.1 are discussed in more detail (and proved)
in the Section 3 of the paper.

A natural question which arises in connection to Theorem 1.1 is the 
following: if the minimum discussed in the theorem is finite, is it also
possible to reach it as $\Phi^{*} ( a,a^{*} )$ for an element $a$ which
is not $R$-diagonal? Up to present we were not able to settle this 
problem. What we can show is its (non-trivial) equivalence to another
problem, also open, of deciding if a certain freeness condition is implied 
by the equality of two free Fisher informations with respect to subalgebras;
see Sections 3.10, 3.11 below.

It is interesting that one can formulate a ``matrix version'' of the 
Theorem 1.1 -- i.e. a version where ``$a$'' becomes a
$d \times d$-matrix over a $W^{*}$-probability space. The possibility of
making such a generalization is created by the following result,
which is of independent interest:

$\ $

{\bf 1.2 Theorem.} Let $\ncps$ be a $W^{*}$-probability space, with 
$\varphi$ faithful trace, and let $d$ be a positive integer. Then:

\vspace{4pt}

$1^{o}$ For every matrix $A = \mataij \in M_{d} ( \A )$ we have:
\begin{equation}
\Phi^{*} ( \ \aijaijs \ ) \ \geq \ d^{3} \Phi^{*} ( A, A^{*} );
\end{equation}
moreover, (1.4) holds with equality if $\{ A, A^{*} \}$ is free from the
subalgebra of ``scalar matrices'' $M_{d} ( \C I ) \subseteq M_{d} ( \A )$ 
(with $I$ = the unit of $\A$).

\vspace{4pt}

$2^{o}$ For every  selfadjoint matrix $B = ( b_{ij} )_{i,j=1}^{d} \in 
M_{d} ( \A )$ we have:
\begin{equation}
\Phi^{*} ( \ \{ b_{ij} \}_{1 \leq i,j \leq d} \ ) \ \geq \
d^{3} \Phi^{*} (B);
\end{equation}
and (1.5) holds with equality if $B$ is free from 
$M_{d} ( \C I ) \subseteq M_{d} ( \A )$.

$\ $

It is easy to see that the freeness conditions appearing in 
Theorem 1.2 can indeed be fulfilled, in the context where the 
$*$-distribution of $A$ (in $1^{o}$) and the distribution of $B$
(in $2^{o}$) are prescribed -- see the discussion preceding 
Proposition 4.1 in Section 4.

The conditions under which equality is reached in (1.4), (1.5) have 
again to do with the more general concept of free Fisher information 
with respect to a subalgebra. For instance, the fact standing behind 
the statement of Theorem $1.2.1^{o}$ is the following: if in addition 
to the family $\aijs$ we also consider a unital $W^{*}$-subalgebra 
$\B \subseteq \A$, then:
\begin{equation}
\Phi^{*} ( \ \aijaijs : \B \ ) \ = \ d^{3} 
\Phi^{*} ( \ \{ A,A^{*} \} : \M_{d} ( \B ) \ );
\end{equation}
in the particular case when $\B = \C I$, this leads to
\[
\Phi^{*} ( \ \aijaijs  \ ) \ = \ d^{3} 
\Phi^{*} ( \ \{ A,A^{*} \} : \M_{d} ( \C I ) \ ) \ \geq \
d^{3} \Phi^{*} ( \ \{ A,A^{*} \} \ ),
\]
which is (1.4) (see Proposition 4.1 below, and the comment following
to it).

By combining the results of the Theorems 1.1 and $1.2.1^{o}$, one obtains
the above mentioned generalization of 1.1:

$\ $

{\bf 1.3 Theorem.} Let $\nu$ and $\mu$ be as in Theorem 1.1, and let $d$ 
be a positive integer. Then
\begin{equation}
\min  \Bigl\{ \ \Phi^{*} ( \ \aijaijs \ ) \ |  \ 
\begin{array}{l} 
A := \mataij \in M_{d} ( \A ) \mbox{ is such }  \\
\mbox{that } A^{*}A \mbox{ has distribution $\nu$ }
\end{array} \ \Bigr\} \ = \ 2 d^{3} \Phi^{*} ( \mu ).
\end{equation}
The minimum is attained if the matrix $A = \mataij$ is an $R$-diagonal
element of $M_{d} ( \A )$, and if $\{ A, A^{*} \}$ is free from the algebra
of scalar matrices $M_{d} ( \C I ) \subseteq M_{d} ( \A )$.

$ \ $

It is easy to see that minimization problems for $\Phi^{*}$ correspond to
maximization problems for the concept of free entropy $\chi^{*}$, which
was also defined (in terms of $\Phi^{*}$) in Voiculescu's work 
\cite{V-P5}. We will conclude the paper by spelling out the maximization
results for $\chi^{*}$ which follow from the theorems presented above.
The counterpart of Theorem 1.3 is:

$\ $

{\bf 1.4 Theorem.} Let $\nu$ and $\mu$ be as in Theorem 1.1, and let $d$ 
be a positive integer. Then
\[
\max  \Bigl\{ \ \chi^{*} ( \ \aijaijs \ ) \ |  \ 
\begin{array}{l} 
A := \mataij \in M_{d} ( \A ) \mbox{ is such }  \\
\mbox{that } A^{*}A \mbox{ has distribution $\nu$ }
\end{array} \ \Bigr\}
\]
\begin{equation}
= \ 2 d^{2} \Bigl( \chi^{*} ( \mu ) - \frac{\log d}{2} \Bigr) .
\end{equation}
The maximum is attained if the matrix $A = \mataij$ is an $R$-diagonal
element of $M_{d} ( \A )$, and if $\{ A, A^{*} \}$ is free from the
algebra of scalar matrices $M_{d} ( \C I ) \subseteq M_{d} ( \A )$.

$\ $

The Theorem 1.4 is obtained from its particular case $d=1$ via a 
maximization result for the free entropy of a family of matrix entries,
which is an analogue of Theorem 1.2:

$\ $

{\bf 1.5 Theorem.} Let $\ncps$ be a $W^{*}$-probability space, with 
$\varphi$ faithful trace, and let $d$ be a positive integer. Then:

\vspace{4pt}

$1^{o}$ For every matrix $A = \mataij \in M_{d} ( \A )$ we have:
\begin{equation}
\chi^{*} ( \ \aijaijs \ ) \ \leq \ d^{2} \Bigl( \chi^{*} 
( A, A^{*} ) - \log d \Bigr) ;
\end{equation}
moreover, (1.9) holds with equality if $\{ A, A^{*} \}$ is free from the
subalgebra of scalar matrices $M_{d} ( \C I ) \subseteq M_{d} ( \A )$.

\vspace{4pt}

$2^{o}$ For every  selfadjoint matrix $B = ( b_{ij} )_{i,j=1}^{d} \in 
M_{d} ( \A )$ we have:
\begin{equation}
\chi^{*} ( \ \{ b_{ij} \}_{1 \leq i,j \leq d} \ ) \ \leq \
d^{2} \Bigl( \chi^{*} (B) - \frac{\log d}{2} \Bigr) ;
\end{equation}
and (1.10) holds with equality if $B$ is free from 
$M_{d} ( \C I ) \subseteq M_{d} ( \A )$.

$\ $

It is tempting to believe that the results obtained about $\chi^{*}$
in this way remain true if ``$\chi^{*}$'' is replaced by ``$\chi$'', the
free entropy defined via approximations with matrices which was 
studied in \cite{V-P2}--\cite{V-P4}. But at the moment it is not proved
(though it might very well be true) that $\chi$ and $\chi^{*}$ coincide;
and consequently, when we replace $\chi^{*}$ by $\chi$ in our 
maximization results, we just obtain some statements for which proofs are
needed. We hope to discuss these statements about $\chi$ (and supply their
proofs) in a future work.

The paper is organized as follows: after reviewing the concept of free
Fisher information in Section 2, we will prove the Theorem 1.1 in 
Section 3, the Theorems 1.2, 1.3 in Section 4, and the Theorems 1.4, 1.5
in Section 5.

$\ $

$\ $

{\bf Acknowledgement:} Part of the research reported in this paper was 
done during a ``Research in Pairs'' program (supported by 
Volkswagen Stiftung) of the Mathematisches Forschungsinstitut 
Oberwolfach, Germany. We would like to acknowledge the excellent work
conditions and very stimulating atmosphere provided by this program.

$\ $

$\ $

$\ $

\setcounter{section}{2}
\setcounter{equation}{0}
{\large\bf 2. Review of the concept of free Fisher information}

$\ $

For general ``free probabilistic'' terminology and basic results,
we refer the reader to the monograph \cite{VDN}.

$\ $

{\bf 2.1 Notations.} Let $\ncps$ be a $W^{*}$-probability space, with 
$\varphi$ a faithful trace.

\vspace{4pt}

$1^{o}$ $\laphi$ will denote the Hilbert space obtained by completing
$\A$ with respect to the norm $|| a ||_{2} := \sqrt{\varphi (a^{*}a)}$,
$a \in \A$.

\vspace{4pt}

$2^{o}$ For $d$ a positive integer, we will denote by $M_{d} ( \A )$ 
the $W^{*}$-algebra of $d\times d$-matrices over $\A$. Also, we will 
denote: $\varphi_{d} \ := \  tr \otimes \varphi : 
M_{d} ( \A ) \rightarrow \C$, where $tr$ is the normalized trace on 
$M_{d} ( \C )$. In other words, $\varphi_{d}$ is 
the faithful trace-state which acts by the formula
\begin{equation}
\varphi_{d} (A) \ = \ \frac{1}{d} \sum_{i=1}^{d} \varphi (a_{ii}), \
\mbox{ for } A = \mataij \in M_{d} ( \A ).
\end{equation}

\vspace{4pt}

$3^{o}$ An immediate consequence of (2.1) is that:
\begin{equation}
|| \ A  \ ||_{\lphid}^{2} \ = \ \frac{1}{d} \sum_{i,j=1}^{d} 
|| \ a_{ij} \ ||_{\lphi}^{2}, \ \ \forall A = \mataij \in M_{d} ( \A ).
\end{equation}
Thus if we fix a pair of indices $k,l \in \{ 1, \ldots , d \}$, then
we get:
\[
|| \ a_{k,l} \ ||_{\lphi} \ \leq \ \sqrt{d}  \ || \ A \ ||_{\lphid},
\ \ \forall A = \mataij \in M_{d} ( \A );
\]
and consequently, the map $A \mapsto a_{k,l}$ extends by continuity to a 
bounded linear map ``$\Entry_{k,l}$'' from $\ldaphi$ to $\laphi$. 
Equation (2.2) can then be extended by continuity to:
\begin{equation}
|| \ X  \ ||_{\lphid}^{2} \ = \ \frac{1}{d} \sum_{i,j=1}^{d} 
|| \ \Entry_{i,j} (X) \ ||_{\lphi}^{2}, \ \ \forall X \in \ldaphi ;
\end{equation}
and by using (2.3) it is readily seen that $X \mapsto
( \Entry_{i,j} (X) )_{i,j=1}^{d}$ is a bijection between $\ldaphi$ and the 
vector space of $d\times d$-matrices over $\laphi$. We will identify in what
follows the vectors in $\ldaphi$ with matrices over $\laphi$, via this 
bijection. It is easily checked that, in this identification, the left and
right actions of $M_{d} ( \A )$ on $\ldaphi$ become ``matrix 
multiplications'' -- e.g. we have that:
\[
\Entry_{k,l} (XA) \ = \ \sum_{m=1}^{d} \Entry_{k,m} (X) \cdot a_{m,l},
\]
for every $X \in \ldaphi$, $A = \mataij$, $1 \leq k,l \leq d$.
The formulas for the entries of $X^*$, and for $\ff_d(X)$, 
$X\in \ldaphi$, are also obtained by continuity in the obvious way.

$\ $

The considerations made in this paper revolve around the notion of free
Fisher information, which was introduced and studied in \cite{V-P1},
\cite{V-P5}. We will next review this notion (Sections 2.2--2.6). A family 
$\ais$ of elements of a $W^{*}$-algebra will be called in what follows 
``selfadjoint'' if there exists an involutive bijection 
$\sigma : I \rightarrow I$ such that $a_{i}^{*} = a_{\sigma (i) }$ for 
every $i \in I$.

$\ $

{\bf 2.2 Definition.} Let $\ncps$ be a $W^{*}$-probability space, with
$\varphi$ faithful trace. Let $\{ a_{i} \}_{i \in I}$ be a selfadjoint
family of elements of $\A$, and let $\B \subseteq \A$ be a unital
$W^{*}$-algebra. 

\vspace{4pt}

$1^{o}$ We say that a family $\xiis$ of vectors in $\laphi$
fulfills the conjugate relations for $\ais$, with respect to $\B$,
if:
\begin{equation}
\varphi ( \xi_{i} b_{0} a_{i_{1}} b_{1} \cdots a_{i_{n}} b_{n} ) \ = \ 
\sum_{m=1}^{n} 
\delta_{i,i_{m}} \varphi ( b_{0} a_{i_{1}} \cdots a_{i_{m-1}} b_{m-1} ) 
\cdot \varphi ( b_{m} a_{i_{m+1}} \cdots a_{i_{n}} b_{n} ), 
\end{equation}
for every $n \geq 0$, $b_{0},b_{1}, \ldots , b_{n} \in \B$ and 
$i, i_{1}, \ldots , i_{n} \in I$.

\vspace{4pt}

$2^{o}$ We say that a family $\xiis$ of vectors in $\laphi$ is a 
conjugate system for $\ais$ with respect to $\B$  if it satisfies 
the Equation (2.4) and if in addition we have that:
\begin{equation}
\xi_{i} \in  \overline{ Alg( { \{ a_{j} \} }_{j \in I} \cup 
\B ) }^{|| \cdot ||_{2}}
\ \subseteq \ L^{2} \ncps , \ \ \forall i \in I .
\end{equation}

$\ $

{\bf 2.3 Remarks.} $1^{o}$ The conjugate relations (2.4) can be viewed 
as a prescription for the inner products in $\laphi$ between a $\xi_{i}$
$(i \in I)$ and a monomial $b_{0} a_{i_{1}} b_{1} \cdots a_{i_{n}} b_{n}$; 
since the monomials of this form linearly span $Alg( \ais \cup \B )$,
it follows that the conjugate system $\xiis$ for $\ais$ with respect to 
$\B$ is unique, if it exists. Note moreover that the existence of the 
conjugate system is equivalent to the existence of any family of vectors 
in $\laphi$ which fulfill the conjugate relations (2.4); indeed, if
$\xiis$ satisfy (2.4) and if we set $\eta_{i}$ to be the projection of 
$\xi_{i}$ onto
$\overline{ Alg( \{ a_{j} \}_{j \in I} \cup  \B ) }^{|| \cdot ||_{2}}$, 
$i \in I$, then $\etais$ will also satisfy (2.4), hence will give the 
conjugate system.

\vspace{4pt}

$2^{o}$ If the family $\ais$ from Definition 2.2 has a conjugate system
$\xiis$ with respect to $\B$, and if $\sigma : I \rightarrow I$ is an 
involution such that $a_{i}^{*} = a_{\sigma (i) }$, $i \in I$,
then we necessarily also have:
\begin{equation}
\xi_{i}^{*} \ = \ \xi_{\sigma (i)}, \ \ i \in I.
\end{equation}
Indeed, it is easy to see (by using the relations $a_{i}^{*} = 
a_{\sigma (i)}$, $i \in I$, and the properties of the trace-state $\varphi$)
that if we set $\eta_{i} = \xi_{\sigma (i)}^{*}$, $i \in I$, then $\etais$
will also fulfill the conjugate relations (2.4); therefore $\eta_{i} =
\xi_{i}$, $i \in I$, by the uniqueness of the conjugate system, and this 
gives (2.6).

$\ $

{\bf 2.4 Definition.} Let $\ncps$ be a $W^{*}$-probability space, with
$\varphi$ faithful trace, let $\{ a_{i} \}_{i \in I}$ be a selfadjoint 
family of elements of $\A$, and let $\B \subseteq \A$ be a unital
$W^{*}$-subalgebra. If $\ais$ has a conjugate system $\xiis$ with respect
to $\B$, then the free Fisher information of $\ais$ with respect to
$\B$ is:
\begin{equation}
\Phi^{*} ( \ \ais \ : \B ) \ := \ \sum_{i \in I} || \xi_{i} ||^{2} .
\end{equation}
If $\ais$ has no conjugate system with respect to $\B$, then one takes
$\Phi^{*} ( \ \ais \ : \B ) \ := \ \infty$.

$\ $

{\bf 2.5 Definition.} Let $\ncps$ be a $W^{*}$-probability space, with 
$\varphi$ faithful trace. If $\ais$ is a selfadjoint family of elements
of $\A$, then we denote:
\begin{equation}
\Phi^{*} ( \ \ais \ ) \ := \ \Phi^{*} ( \ \ais \ : \C I ).
\end{equation}
$\Phi^{*} ( \ \ais \ )$ will be simply called ``the free Fisher
information'' of $\ais$. Also, if $\xiis$ fulfills the conjugate relations
(respectively is a conjugate system) for $\ais$ with respect to $\C I$,
we will generally omit ``with respect to $\C I$'' from the formulation.

$\ $

{\bf 2.6 Remarks.} $1^{o}$ Let $\ncps$, $\ais$ and $\B$ be as in the 
Definition 2.4. If a family $\xiis$ in $\laphi$ fulfills the conjugate
relations for $\ais$ with respect to $\B$, but does not necessarily
satisfy (2.5), then we still know that:
\begin{equation}
\Phi^{*} ( \ \ais \ : \B ) \ \leq \ \sum_{i \in I} || \xi_{i} ||^{2}.
\end{equation}
This is a direct consequence of the statement concluding the Remark
$2.3.1^{o}$.

\vspace{4pt}

$2^{o}$ Let $\ncps$ and $\ais$ be as above, and let 
$\B_{1}, \B_{2}$ be $W^{*}$-subalgebras of $\A$ such that
$I \in \B_{1} \subseteq \B_{2}$. Then
\begin{equation}
\Phi^{*} ( \ \ais \ : \B_{1} ) \ \leq \ \Phi^{*} ( \ \ais \ : \B_{2} ).
\end{equation}
Indeed, if $\Phi^{*} ( \ \ais \ : \B_{2} ) < \infty$, then the conjugate
system for $\ais$ with respect to $\B_{2}$ will fulfill the conjugate 
relations with respect to $\B_{1}$; hence (2.10) follows from (2.9).

\vspace{4pt}

$3^{o}$ In the particular case of $2^{o}$ when $\B_{1} = \C I$, we obtain 
the inequality:
\begin{equation}
\Phi^{*} ( \ \ais \ ) \ \leq \ \Phi^{*} ( \ \ais \ : \B ),
\end{equation}
for every unital $W^{*}$-subalgebra $\B$ of $\A$. It is important to 
record here that, as proved in \cite{V-P5} Proposition 3.6, (2.11) holds
with equality whenever $\ais$ is free from $\B$.

$\ $

The problems discussed in the present paper are formulated only in terms
of the free information $\Phi^{*} ( \ \ais \ )$ (with respect to the 
scalars). But however, considerations involving free information with 
respect to non-trivial subalgebras appear naturally in the solutions.
Moreover, in Section 3 we will arrive to use a version of 
$\Phi^{*} ( \ \bullet \ : \B )$ where (in addition to $\B$ itself) one
also considers a completely positive map $\eta : \B \rightarrow \B$.
This version of $\Phi^{*}$ was introduced in \cite{Sh2}, and is defined
as follows.

$\ $

{\bf 2.7 Definition.} Let $\ncps$ be a $W^{*}$-probability space, with
$\varphi$ faithful trace, let $x = x^{*}$ be in $\A$, let
$\B \subseteq \A$ be a unital $W^{*}$-subalgebra, and let
$\eta : \B \rightarrow \B$ be a completely positive map.

\vspace{4pt}

$1^{o}$ We say that a vector $\xi \in \laphi$ fulfills the conjugate 
relations for $x$, with respect to $\B$ and $\eta$, if:
\begin{equation}
\varphi ( \xi b_{0} x b_{1} \cdots x b_{n} ) \ = \ 
\sum_{m=1}^{n} 
\varphi ( \ \eta ( E_{\B} ( b_{0} x \cdots x b_{m-1} )) \cdot 
b_{m} x \cdots x b_{n} ), 
\end{equation}
for every $n \geq 0$ and every $b_{0},b_{1}, \ldots , b_{n} \in \B$,
and where $E_{\B}$ denotes the unique trace-preserving conditional 
expectation from $\A$ onto $\B$.

\vspace{4pt}

$2^{o}$ The vector $\xi \in \laphi$ is called a conjugate for $x$, with 
respect to $\B$ and $\eta$, if it satisfies (2.12) and if in addition:
\begin{equation}
\xi \in  \overline{ Alg( \{ x \} \cup \B ) }^{|| \cdot ||_{2}}.
\end{equation}

\vspace{4pt}

$3^{o}$ The free Fisher information of $x$ with respect to $\B$ and 
$\eta$ is defined to be:
\begin{equation}
\Phi^{*} ( \ x \ : \B , \eta ) \ := \ || \xi ||^{2} ,
\end{equation}
if $x$ has a conjugate vector $\xi$ with respect to $\B$ and $\eta$, 
and $\Phi^{*} ( \ x \ : \B , \eta ) \ := \ \infty$ otherwise.

$\ $

{\bf 2.8 Remarks.} $1^{o}$ Exactly as in Remark $2.3.1^{o}$, one sees
that the conjugate vector with respect to $\B$ and $\eta$ is unique, if
it exists. (This ensures that the definition of 
$\Phi^{*} ( \ x \ : \B , \eta )$ in (2.14) makes sense.)

\vspace{4pt}

$2^{o}$ In the particular case when the completely positive map
$\eta : \B \rightarrow \B$ is $\eta (b) := \varphi (b) I$, $b \in \B$,
one obtains $\Phi^{*} ( \ x \ : \B , \eta ) = \Phi^{*} ( \ x \ : \B )$,
because (2.12) reduces to (2.4).

\vspace{4pt}

$3^{o}$ It is easy to see (exactly as in the Remark $2.6.2^{o}$) that
one has the inequality:
\begin{equation}
\Phi^{*} ( \ x \ : \B_{1} , \eta_{1} ) \ \leq \
\Phi^{*} ( \ x \ : \B_{2} , \eta_{2} )
\end{equation}
whenever $\B_{1} \subseteq \B_{2}$ and $\eta_{1} , \eta_{2}$ are related
by:
\begin{equation}
\eta_{2} (b) \ = \ \eta_{1} ( E_{\B_{1}} (b) ), \ \ \forall b \in \B_{2}.
\end{equation}
It is again important to record that, as proved in \cite{Sh2}
Proposition 3.8, (2.15) holds with equality whenever
$Alg( \{ x \} \cup \B_{1} )$ is free from $\B_{2}$, with amalgamation over
$\B_{1}$.

$\ $

{\bf 2.9 Remark.} Let $( \A_{1} , \varphi_{1} )$ and
$( \A_{2} , \varphi_{2} )$ be $W^{*}$-probability spaces, with
$\varphi_{1} , \varphi_{2}$ faithful traces, and let
$x_{1} = x_{1}^{*}\in A_{1}$, $x_{2} = x_{2}^{*} \in \A_{2}$ be elements 
with identical distributions (i.e., 
$\varphi_{1} ( x_{1}^{n} ) = \varphi_{2} ( x_{2}^{n} ), \forall n \geq 0$). 
Then we must also have that $\Phi^{*} ( x_{1} ) = \Phi^{*} ( x_{2} )$. 
Indeed, the coincidence of distributions has as consequence that there
exists a unitary operator $U: { \overline{ Alg( I,x_{1} ) } }^{|| \cdot 
||_{2}} \rightarrow { \overline{ Alg(I,x_{2}) } }^{|| \cdot ||_{2}}$,
determined by the relation $U ( x_{1}^{n} ) \ = \  x_{2}^{n}$, $n \geq 0$. 
It is immediate that $U$ sends a conjugate for $x_{1}$ into a 
conjugate for $x_{2}$, and this in turn implies the equality of free
Fisher informations.

In particular, if $\mu$ is a probability measure with compact support on
${\bf R}$, it makes sense to use the notation
\begin{equation}
\Phi^{*} ( \mu ) \ := \ \Phi^{*} (x),
\end{equation}
where $x$ is an arbitrary selfadjoint random variable (in some
$W^{*}$-probability space $\ncps$, with $\varphi$ faithful trace) such that
the distribution of $x$ is $\mu$. A detailed discussion about 
$\Phi^{*} ( \mu )$ is made in \cite{V-P1} (see also Section 2 of \cite{V-P5}); 
it is in particular shown there that if $\mu$ is absolutely continuous with 
respect to the Lebesgue measure, and has density $\rho$, then: 
\begin{equation}
\Phi^{*} ( \mu ) \ = \ \frac{2}{3} \cdot 
\int_{- \infty}^{\infty} \rho (t)^{3} \ dt.
\end{equation}

$\ $

$\ $

\setcounter{section}{3}
\setcounter{equation}{0}
{\large\bf 3. Minimization of $\Phi^{*} (a,a^{*})$, when the
distribution of $a^{*}a$ is prescribed} 

$\ $

Let $\nu$ be a probability measure with compact support on $[0, \infty )$. 
We consider the minimization problem stated in (1.1) of the Introduction,
i.e. the problem of determining:
\[
\inf \{ \Phi^{*} ( a, a^{*} ) \ | \ a^{*}a \mbox{ has distribution }
\nu \}
\]
where $a \in \A$ and $\ncps$ is a $W^{*}$-probability space, with 
$\varphi$ faithful trace.

In the considerations related to this problem, it is convenient 
to use the following symmetric measure associated to $\nu$.

$\ $

{\bf 3.1 Definition.} For $\nu$ as above, we will call
``symmetric square root of $\nu$'' the unique probability measure 
$\mu$ on $\R$ which is symmetric (i.e. $\mu (S) = \mu (-S)$ for every 
Borel set $S$), and has the property that 
$\mu (S)$ = $\nu ( \ \{ s^{2} \ | \ s \in S \} \ )$, for every Borel set
$S$ such that $S = -S$.

$\ $

In terms of random variables, the connection between $\nu$ and its 
symmetric square root $\mu$ is expressed as follows: a selfadjoint element 
$x$ in a $W^{*}$-probability space $\ncps$ has distribution $\mu$ if and
only if $x$ is even (i.e. $\varphi ( x^{n} ) = 0$ for $n$ odd), and $x^{2}$ 
has distribution $\nu$.

$\ $

{\bf 3.2 Theorem.} Let $\nu$ be a probability measure with compact support 
on $[0, \infty )$, and let $\mu$ be the 
symmetric square root of $\nu$. Let $\ncps$ be a $W^{*}$-probability 
space, with $\varphi$ faithful trace, and let $a \in \A$ be such that 
$a^{*}a$ has distribution $\nu$. Then:
\begin{equation}
\Phi^{*} ( a, a^{*} ) \ \geq \ 2 \Phi^{*} ( \mu ).
\end{equation}
Moreover, (3.1) holds with equality if $a$ is of the form $a=up$, where
$u \in \A$ is a unitary with Haar distribution (i.e. $\varphi (u^{n}) =0$
for all $n \in \Z \setminus \{ 0 \}$), $p=p^{*}$ has distribution $\mu$,
and $\{ p \}$ is free from $\{ u,u^{*} \}$.

Thus the infimum considered in (1.1) of the Introduction is equal to 
$2 \Phi^{*} ( \mu )$.

$ \ $

{\bf 3.3 Remarks.} $1^{o}$ If $u$ is a unitary with Haar distribution,
$p=p^{*}$, and $\{ p \}$ is free from $\{ u,u^{*} \}$, then the element
$a=up$ is said to be $R$-diagonal (\cite{NS}). For such an element, the 
$*$-distribution of $a$ is completely determined by 
the distribution of $p^{2}$ (\cite{NS}, Corollary 1.8).
This implies that, as far as $*$-distributions are 
concerned, there is a unique $R$-diagonal element $a$ such that the 
distribution of $a^{*}a$ is a given probability measure $\nu$.

Let us hence notice that in the phrase following Equation (3.1) (in the
statement of Theorem 3.2) we could replace ``$p$ has distribution $\mu$''
with the apparently more general condition ``$p^{2}$ has distribution 
$\nu$''. But this wouldn't actually change the $*$-distribution of $a$
-- we would still have to do with the same $R$-diagonal element.

We were in fact unable to determine if the $R$-diagonal $*$-distribution
is the unique one which achieves the minimization of $\Phi^{*} (a,a^{*})$
considered in (1.1). (See also the Sections 3.10, 3.11 below.)

\vspace{4pt}

$2^{o}$ The statement of Theorem 3.2 contains the one of Theorem 1.1, 
with the exception of the formula (1.3). The latter formula follows from
Equation (2.18) of Remark 2.9, combined with the simple observation that
$\mu$ is absolutely continuous if and only if $\nu$ is so, in which case
the densities $\sigma$ of $\mu$ and $\rho$ of $\nu$ are connected by the
relation $\sigma (t) = |t| \rho ( t^{2} )$, $t \in \R$.

$\ $

Our goal in this section is thus to prove Theorem 3.2. Let us set the 
following:

$\ $

{\bf 3.4 Notations.} $\nu$, $\mu$, $\ncps$, $a \in \A$ are fixed from now
on, until the end of the section, and are as in the statement of Theorem 
3.2. We will consider the space $\ncpst$ of $2 \times 2$ matrices over 
$\ncps$ (as in Notations $2.1.2^{o}$), and we will give a special attention 
to the selfadjoint matrix:
\begin{equation}
A \  := \left(  \begin{array}{cc}
0  &   a     \\
a^*  &  0
\end{array}   \right)   \  \in \ M_{2}  ( \A ).
\end{equation}

For $i,j \in \{ 1,2 \}$ we will denote by $V_{ij}$ the matrix in 
$M_{2} ( \A )$ which has the $(i,j)$-entry equal to the unit of $\A$,
and the other entries equal to $0$. Then:
\[
\mbox{span} \{ V_{11}, V_{12}, V_{21}, V_{22} \} \ = \ 
M_{2} ( \C I ) \subseteq M_{2} ( \A );
\]
we will also denote:
\[
\D \ := \ \mbox{span} \{ V_{11}, V_{22} \}
\]
(the 2-dimensional $*$-subalgebra of $M_{2} ( \A )$ consisting of scalar 
diagonal matrices). We will denote by $E_{\M}$ and $E_{\D}$ the unique 
trace-preserving conditional expectations from $M_{2} ( \A )$ onto 
$M_{2} ( \C I )$ and $\D$, respectively. For
$B = ( b_{ij} )_{i,j=1}^{2} \in M_{2} ( \A )$ we have:
\begin{equation}
\expb (B) \ = \
\left(  \begin{array}{cc}
\varphi ( b_{11} ) I  &  \varphi ( b_{12} ) I  \\ 
\varphi ( b_{21} ) I  &  \varphi ( b_{22} ) I 
\end{array}  \right) , \ \ 
\expd (B) \ = \ 
\left(  \begin{array}{cc}
\varphi ( b_{11} ) I  &  0  \\ 
 0                    &  \varphi ( b_{22} ) I
\end{array}  \right) . 
\end{equation}

$\ $

{\bf 3.5 Remark.} Since $A$ of Equation (3.2) has:
$$
A^2 \  = \left(  \begin{array}{cc}
aa^*  &   0     \\
0 &  a^*a
\end{array}   \right),
$$
while on the other hand the odd powers of $A$ have 0's on the main
diagonal, it is immediate that $A$ is even and that $A^2$ has
distribution $\nu$. Therefore $A$ itself has distribution $\mu$.

$\ $

{\bf 3.6 Proposition.} Let $\eta : \MM \rightarrow \MM$ be the 
completely positive map defined by:
\begin{equation}
\eta \ \Bigl( \ \left(   \begin{array}{cc}
x_{11}  &  x_{12}  \\ x_{21}  &  x_{22}  
\end{array}  \right) \ \Bigr) \ :=
\left(   \begin{array}{cc}
x_{22}  &  0  \\ 0  &  x_{11}  
\end{array}  \right) .
\end{equation}
Then we have:
\begin{equation}
\Phi^{*} (a,a^{*}) \ = \ 2 \Phi^{*} (A: \MM , \eta ).
\end{equation}

$\ $

{\bf Proof.} We first consider the situation when
$\Phi^* (a,a^* ) < \infty$. In this case there exists $\xi \in
\overline{Alg( I,a,a^{*} )}^{|| \cdot ||_{2}}$ such that 
$\{ \xi , \xi^* \}$ forms a conjugate system for $\{ a,a^* \}$.
We define:
\begin{equation}
X \ := \left(  \begin{array}{cc}
0  &   \xi^*     \\
\xi  &  0
\end{array}   \right)   \  \in \ \laphitwo
\end{equation}
where the identification between vectors in $\laphitwo$ and matrices
over $\laphi$ is as discussed in the Notations $2.1.3^{o}$. We will
show that $X$ is a conjugate for $A$, with respect to $\MM$ and $\eta$.
Proving this claim consists in verifying that:

(a) the relation
\begin{equation}
\varphi_{2} ( XB_{0}AB_{1} \cdots AB_{n} ) \ = \
\sum_{m=1}^{n}
\varphi_{2} ( \ \eta ( \expb ( B_{0}A \cdots AB_{m-1} ) \cdot 
B_{m}A \cdots AB_{n} ) 
\end{equation}
holds for every $n \geq 0$ and every $B_{0},B_{1}, \ldots , B_{n} \in \MM$;
and that:

\vspace{6pt}

(b) $X \in \overline{Alg( \{ A \} \cup \MM )}^{|| \cdot ||_{2}}
\subseteq \laphitwo$.

\vspace{6pt}

\noindent
Note that once (a) and (b) will be proved, we will have the equality
\[
\Phi^*(A: \MM , \eta ) =  \Vert X\Vert^2_{\lphitwo}
\stackrel{(2.3)}{=} \frac 12 (\Vert\xi\Vert^2+\Vert\xi^*\Vert^2)
= \frac 12 \Phi^*(a,a^*),
\]
which is exactly (3.5) (under the hypothesis $\Phi^{*} (a,a^{*})
< \infty$).

\vspace{10pt}

{\em Proof of (a).} Both sides of (3.7) depend multilinearly on 
$B_{0}, B_{1}, \ldots , B_{n}$; we can therefore assume without loss of
generality that $B_{m} = V_{i_{m}j_{m}}, \ 0 \leq m \leq n$, for some
$i_{0},j_{0}, \ldots , i_{n}, j_{n}$ $\in$ $\{ 1,2 \}$.

By using the trace property of $\varphi_{2}$ we can write the left-hand
side of (3.7) as:
\begin{equation}
\varphi_{2} ( \ V_{j_{n}j_{n}} X V_{i_{0}j_{0}} A 
V_{i_{1}j_{1}}  \cdots A V_{i_{n}j_{n}} \ ).
\end{equation}
Only the $(j_{n},j_{n})$-entry of the matrix product appearing in (3.8)
is not 0; this entry equals:
\[
(X)_{j_{n}i_{0}} (A)_{j_{0}i_{1}} \cdots (A)_{j_{n-1}i_{n}},
\]
where $(A)_{ij}, (X)_{ij}$ stand for the $(i,j)$-entry of $A$ and $X$,
respectively. Thus the quantity in (3.8) equals:
\begin{equation}
\frac{1}{2} \varphi ( \ 
(X)_{j_{n}i_{0}} (A)_{j_{0}i_{1}} \cdots (A)_{j_{n-1}i_{n}} \ ).
\end{equation}
But we know that $(A)_{ij} = 0 = (X)_{ij}$ if $i=j$; so if we make the 
convention to denote $\overline{i} :=3-i$ ( = the number in $\{ 1,2 \}$ 
which is not $i$), for $i \in \{ 1,2 \}$, then (3.9) becomes:
\begin{equation}
\frac{1}{2}  \delta_{j_{0}{\barl{i}}_{1}}
\delta_{j_{1}{\barl{i}}_{2}} \cdots \delta_{j_{n-1}{\barl{i}}_{n}}
\delta_{j_{n}{\barl{i}}_{0}} 
\cdot \varphi ( \ (X)_{ {\barl{i}}_{0} i_{0}} 
(A)_{ {\barl{i}}_{1} i_{1}} \cdots  (A)_{ {\barl{i}}_{n} i_{n}}  \ ).
\end{equation}
In (3.10), $(X)_{ {\barl{i}}_{0} i_{0}}$ is $\xi$ or $\xi^{*}$, while
every $(A)_{ {\barl{i}}_{m} i_{m}}$ is either $a$ or $a^{*}$. So the 
conjugate relations for $\{ a,a^{*} \}$ can  be used, to obtain that the 
quantity in (3.10) equals:
\[
\frac{1}{2}  \delta_{j_{0}{\barl{i}}_{1}}
\delta_{j_{1}{\barl{i}}_{2}} \cdots \delta_{j_{n-1}{\barl{i}}_{n}}
\delta_{j_{n}{\barl{i}}_{0}} \cdot
\]
\begin{equation}
\cdot \sum_{m=1}^{n} \delta_{ {\barl{i}}_{0} i_{m} }
\varphi ( \ (A)_{ {\barl{i}}_{1} i_{1}} \cdots  
(A)_{ {\barl{i}}_{m-1} i_{m-1}}  \ ) \cdot
\varphi ( \ (A)_{ {\barl{i}}_{m+1} i_{m+1}} \cdots  
(A)_{ {\barl{i}}_{n} i_{n}}  \ ).
\end{equation}

We now turn to the right-hand side of (3.7). By using the formulas for 
$\eta$ and $\expb$ (as in Equations (3.4) and (3.3)) we first see that:
\[
\eta ( \expb ( V_{i_{0}j_{0}} A \cdots A V_{i_{m-1}j_{m-1}} )) \ = \
\delta_{i_{0}j_{m-1}} \cdot \varphi ( \ (A)_{j_{0}i_{1}} \cdots
(A)_{j_{m-2}i_{m-1}} \ ) \cdot V_{ {\barl{i}}_{0} {\barl{i}}_{0} }, 
\]
for $1 \leq m \leq n$. By replacing this into the right-hand side of (3.7),
we obtain the expression:
\begin{equation}
\sum_{m=1}^{n} \delta_{i_{0}j_{m-1}} \cdot
\varphi ( \ (A)_{j_{0}i_{1}} \cdots (A)_{j_{m-2}i_{m-1}} \ ) \cdot 
\varphi_{2} ( \ V_{ {\barl{i}}_{0} {\barl{i}}_{0} } V_{i_{m}j_{m}} A
\cdots A V_{i_{n}j_{n}} \ ).
\end{equation}
But then a calculation very similar to the ones shown above gives us that
the summation in (3.12) coincides, term by term, with the one in (3.11).

\vspace{10pt}

{\em Proof of (b).} We have
\begin{equation}
\left(  \begin{array}{cc}
p(a,a^{*})  &  0   \\   0   &   0   \end{array}  \right) \
\in \ Alg( \{ A \} \cup \MM ),
\end{equation}
whenever $p$ is a non-commutative polynomial of two variables. (Indeed,
it clearly suffices to check the cases $p( a,a^{*} ) =I$ and 
$p( a,a^{*} ) = a$, when the matrix in (3.13) becomes $V_{11}$ and 
respectively $AV_{21}$.) From (3.13) and the fact that
$\xi , \xi^{*} \in Alg( I,a,a^{*} )$ we infer:
\begin{equation}
\left(  \begin{array}{cc}
\xi  &  0   \\   0   &   0   \end{array}  \right) , \
\left(  \begin{array}{cc}
\xi^{*}  &  0   \\   0   &   0   \end{array}  \right)  \
\in \ \overline{Alg( \{ A \} \cup \MM )}^{|| \cdot ||_{2}} .
\end{equation}
But $\overline{Alg( \{ A \} \cup \MM )}^{|| \cdot ||_{2}}$
is invariant under the left/right action of elements from
$\MM$; so (3.14) implies that:
\[
X \ = \ V_{21}
\left(  \begin{array}{cc}
\xi  &  0   \\   0   &   0   \end{array}  \right) \ + \
\left(  \begin{array}{cc}
\xi^{*}  &  0   \\   0   &   0   \end{array}  \right) V_{12} \
\in \ \overline{Alg( \{ A \} \cup \MM )}^{|| \cdot ||_{2}} ,
\]
as desired.

\vspace{10pt}

Hence (3.5) is now proved in the case when $\Phi^{*} (a,a^{*}) < \infty$.
It remains to show that $\Phi^{*} (a,a^{*}) = \infty$ $\Rightarrow$
$\Phi^{*} (A: \MM , \eta ) = \infty$; or equivalently that 
$\Phi^{*} (A: \MM , \eta ) < \infty$ $\Rightarrow$
$\Phi^{*} (a,a^{*}) < \infty$.

If $\Phi^{*} (A: \MM , \eta ) < \infty$, then there exists $X \in \laphitwo$
which satisfies the conjugate relations for $A$, with respect to $\MM$ and 
$\eta$. We identify $X$ with a $2 \times 2$ matrix over $\laphi$, and 
denote its (2,1)-entry by $\xi$; we will show that $\{ \xi , \xi^{*} \}$
satisfy the conjugate relations with respect to $\{ a,a^{*} \}$ (this 
will entail, as noticed in Remark $2.6.1^{o}$, that
$\Phi^{*} (a,a^{*} ) < \infty$).

It is in fact sufficient to prove that:
\begin{equation}
\varphi ( \xi a_{i_{1}} \cdots a_{i_{n}} ) \ = \ 
\sum_{m=1}^{n} \delta_{i_{m},1} 
\varphi ( a_{i_{1}} \cdots a_{i_{m-1}} ) 
\varphi ( a_{i_{m+1}} \cdots a_{i_{n}} ), 
\end{equation}
for every $n \geq 1$ and every $i_{1}, \ldots , i_{n} \in \{ 1,2 \}$, 
where we denoted $a_{1} :=a$, $a_{2} := a^{*}$. Indeed, the symmetric
relation:
\begin{equation}
\varphi ( \xi^{*} a_{i_{1}} \cdots a_{i_{n}} ) \ = \ 
\sum_{m=1}^{n} \delta_{i_{m},2} 
\varphi ( a_{i_{1}} \cdots a_{i_{m-1}} ) 
\varphi ( a_{i_{m+1}} \cdots a_{i_{n}} ) 
\end{equation}
follows from (3.15) by taking an adjoint and doing a circular permutation
under $\varphi$. We also have $\varphi ( \xi ) = 2 \varphi_{2} (XV_{12} )
=0$, by the conjugate relations satisfied by $X$, and 
$\varphi ( \xi^{*} ) = \overline{\varphi ( \xi )} = 0$. Added to 
(3.15-16), this exhausts the list of conjugate relations for $a,a^{*}$.

In order to verify (3.15), we adopt again the conventions of notation used
in the ``Proof of (a)'' above, and we write:
\[
\varphi ( \xi a_{i_{1}} a_{i_{2}} \cdots a_{i_{n}} )
\]
\[
\ = \ \varphi (  \ (X)_{21} (A)_{ i_{1} {\barl{i}}_{1} } 
(A)_{i_{2} {\barl{i}}_{2}} \cdots (A)_{ i_{n} {\barl{i}}_{n} }  \ )
\]
\[
\ = \ 2 \varphi_{2} (  X V_{1i_{1}} A V_{ {\barl{i}}_{1} i_{2} } A \cdots 
V_{ {\barl{i}}_{n-1} i_{n} } A V_{ {\barl{i}}_{n} 2 } )
\]
\[
= \ 2 \sum_{m=1}^{n} \varphi_{2} ( \ \eta ( \expb (
V_{1i_{1}} A V_{ {\barl{i}}_{1} i_{2} } A \cdots 
V_{ {\barl{i}}_{m-1} i_{m} } )) \cdot 
V_{ {\barl{i}}_{m} i_{m+1} } A \cdots
V_{ {\barl{i}}_{n-1} i_{n} } A V_{ {\barl{i}}_{n} 2} \ )
\]
(by the conjugate relations for $A$, with respect to $\MM$ and $\eta$)
\[
= \ 2 \sum_{m=1}^{n} \varphi_{2} ( \ \delta_{i_{m},1} \varphi (
(A)_{i_{1} {\barl{i}}_{1}} \cdots (A)_{i_{m-1} {\barl{i}}_{m-1}} ) V_{22} 
\cdot V_{ {\barl{i}}_{m} i_{m+1} } A \cdots
V_{ {\barl{i}}_{n-1} i_{n} } A V_{ {\barl{i}}_{n} 2} \ )
\]
(by writing explicitly how $\eta \circ \expb$ works)
\[
= \ 2 \sum_{m=1}^{n} \delta_{i_{m},1} \varphi (
(A)_{i_{1} {\barl{i}}_{1}} \cdots (A)_{i_{m-1} {\barl{i}}_{m-1}} ) \cdot 
\frac{1}{2} \delta_{2, {\barl{i}}_{m}} \varphi (
(A)_{i_{m+1} {\barl{i}}_{m+1}} \cdots (A)_{i_{n} {\barl{i}}_{n}} ) 
\]
\[
= \ \sum_{m=1}^{n} \delta_{i_{m},1} 
\varphi ( a_{i_{1}} \cdots a_{i_{m-1}} ) 
\varphi ( a_{i_{m+1}} \cdots a_{i_{n}} ). {\bf QED}
\]

$\ $

\vspace{14pt}

{\bf 3.7 Proposition.} Let $\etaz : \D \rightarrow \D$ be the
$*$-automorphism defined by:
\begin{equation}
\etaz \ \Bigl( \ \left(   \begin{array}{cc}
x_{11}  &  0  \\ 0  &  x_{22}  
\end{array}  \right) \ \Bigr) \ := \
\left(   \begin{array}{cc}
x_{22}  &  0  \\ 0  &  x_{11}  
\end{array}  \right) .
\end{equation}
Then we have:
\begin{equation}
\Phi^{*} (A: \MM , \eta ) \ \geq \ 
\Phi^{*} (A: \D , \etaz ) \ = \
\Phi^{*} (A) \ = \ \Phi^{*} ( \mu ).
\end{equation}

$\ $

{\bf Proof.} It is immediate that $\eta (B) = \etaz ( \expd (B))$, for 
every $B \in \MM$; thus the inequality appearing in (3.18) is implied by
(2.15) of Remark $2.8.3^{o}$. On the other hand, the equality 
$\Phi^{*} (A) = \Phi^{*} ( \mu )$ holds just because $A$ has distribution
$\mu$ (Remark 3.5). Our main concern in this proof is to show that
$\Phi^{*} (A: \D , \etaz ) = \Phi^{*} (A)$.

For proving $\Phi^{*} (A: \D , \etaz ) \geq \Phi^{*} (A)$, we assume
the existence of a conjugate $X$ for $A$, with respect to $\D$ and $\etaz$,
and we show that that $X$ fulfills the conjugate relations for $A$ with
respect to the scalars. The assumption on $X$ is that:
\begin{equation}
\varphi_{2} (XD_{0}AD_{1} \cdots AD_{n} ) \ = \
\sum_{m=1}^{n} \varphi_{2} ( \ \etaz ( \expd (D_{0}A \cdots AD_{m-1} ))
\cdot D_{m}A \cdots AD_{n} \ ),
\end{equation}
for every $n \geq 0$ and every $D_{0},D_{1}, \ldots , D_{n} \in \D$.
By setting in (3.19) $D_{0} = D_{1} = \cdots D_{n} = I_{2}$ (the unit  
of $M_{2} ( \A )$ ), we get:
\begin{equation}
\varphi_{2} (XA^{n} ) \ = \ \sum_{m=1}^{n} \varphi_{2}
( \ \etaz ( \expd (A^{m-1} )) \cdot A^{n-m} \ ), \ \ n \geq 0.
\end{equation}
It is however immediately checked that:
\[
\etaz ( \expd (A^{k} )) \ = \ \varphi_{2} ( A^{k} ) I_{2}, \ \ k \geq 0;
\]
hence (3.20) comes to:
\[
\varphi_{2} (XA^{n} ) \ = \ \sum_{m=1}^{n} \varphi_{2} ( A^{m-1} ) \cdot
\varphi_{2} ( A^{n-m} ), \ \ n \geq 0,
\]
which says exactly that $X$ fulfills the conjugate relations for $A$ with
respect to the scalars.

We now go to the proof of the opposite inequality,
$\Phi^{*} (A: \D , \etaz ) \leq \Phi^{*} (A)$. The method is the same as 
above (although the calculations will be more complicated): we assume that
$A$ has a conjugate vector $X \in \laphitwo$, with respect to the scalars,
and we will show that $X$ also fulfills the conjugate relations for $A$
with respect to $\D$ and $\etaz$. We identify the vector $X$ with a 
matrix over $\laphi$ (as in $2.1.3^{o}$):
\begin{equation}
X \  = \left(  \begin{array}{cc}
\xi_{11}  &   \xi_{12}     \\
\xi_{21}  &   \xi_{22}
\end{array}   \right) ,  \qquad\mbox{with $\xi_{ij}\in\laphi$.}
\end{equation}
Note that 
\begin{eqnarray*}
\mbox{$A = A^*$ in $M_{2} ( \A )$} & \Rightarrow &
\mbox{$X = X^*$ in $\laphitwo$}
\qquad\mbox{(by Remark $2.3.2^{o}$)}
\\
& \Rightarrow &
\mbox{$\xi_{12} = \xi_{21}^{*}$ in $\laphi$.}
\end{eqnarray*}

Before doing anything else, let us show that in (3.21) we have 
$\xi_{11} = \xi_{22} = 0$. To this end we will use ``the even half'' of
the conjugate relations fulfilled by $X$:
\[
\ff(XA^{2k})=\sum_{l=1}^{2k} \ff(A^{l-1})\cdot\ff(A^{2k-l}) , 
\ \ k \geq 0.
\]
Every term in the latter sum is 0, because one of $A^{l-1}$ and
$A^{2k-l}$ must always have vanishing diagonal entries. So we get
$\ff(XA^{2k})=0$, hence $X\perp A^{2k}$ in $\laphitwo$, for every
$k \geq 0$. Since on the other hand the definition of the conjugate vector 
contains the fact that
$$
X\in \overline{\mbox{span}}^{\Vert\cdot\Vert_2}\{
A^n\mid n\geq 0\}\subseteq \laphitwo,
$$
and since (obviously) $A^{n} \perp A^{m}$ when $n,m$ have different 
parities, we infer that actually:
\begin{equation}
X\in \overline{\mbox{span}}^{\Vert\cdot\Vert_2}\{
A^{2k+1}\mid k\geq 0\}=
\overline{\mbox{span}}^{\Vert\cdot\Vert_2}\{
\left(
\begin{array}{cc}
0 & a(a^*a)^k \\ a^*(aa^*)^k & 0 
\end{array} \right)\mid k \geq 0 \} .
\end{equation}
From the discussion in $2.1.3^{o}$ it is clear that convergence
in $\laphitwo$ implies ``entry-wise convergence'' in $\laphi$.
Therefore (3.22) has as consequence that $\xi_{11} = \xi_{22} = 0$, as
desired, and we can write:
\begin{equation}
X=\left(
\begin{array}{cc}
0 & \xi^* \\ \xi & 0 
\end{array} \right),
\end{equation}
where $\xi := \xi_{21}$ of (3.21).

Besides (3.23), there is another consequence of (3.22) which will be
used in the sequel, namely that:
\begin{equation}
\varphi ( \xi a (a^{*}a)^{m} ) \ = \
\varphi ( \xi^{*} a^{*} (aa^{*})^{m} ), \ \ \forall m \geq 0.
\end{equation}
Indeed, for every given $m \geq 0$, (3.22) implies:
\[
\left(  \begin{array}{cc}
\xi^{*} a^{*} (aa^{*})^{m}   &    0   \\
0       &    \xi a (a^{*}a)^{m}   
\end{array}   \right)  \ = \ XA^{2m+1}
\ \in \ {\overline{\mbox{span}}}^{|| \cdot ||_{2}} \{ A^{2k} \ | \
k \geq m+1 \}
\]
\[
= \ {\overline{\mbox{span}}}^{|| \cdot ||_{2}} \{  \ 
\left(  \begin{array}{cc}
(aa^{*})^{k}   &    0   \\
0       &    (a^{*}a)^{k}   
\end{array}   \right) \ | \ k \geq m+1 \ \} . 
\]
Then the fact that $\varphi ( \ (aa^{*})^{k} \ ) =
\varphi ( \ (a^{*}a)^{k} \ )$, $k \geq m+1$, can be passed through the 
closed linear span to yield (3.24).

Now, recall that our goal is to prove that $X$ fulfills the conjugate 
relations for $A$, with respect to $\D$ and $\etaz$; these relations are
exactly as described in (3.19). Since $\D$ = span$\{ V_{11}, V_{22} \}$, 
it actually suffices to check that:
\begin{equation}
\varphi_{2} ( X V_{ i_{0}i_{0} } A V_{ i_{1}i_{1} } \cdots
A V_{ i_{n}i_{n} } ) \ = \ \sum_{m=1}^{n} \varphi_{2} 
( \ \etaz ( \expd ( V_{ i_{0}i_{0} } A \cdots A V_{ i_{m-1}i_{m-1} } ))
\cdot V_{ i_{m}i_{m} } A \cdots A V_{ i_{n}i_{n} } \ ),
\end{equation}
for every $n \geq 0$ and every $i_{0},i_{1}, \ldots , i_{n} \in \{ 1,2 \}$.

The verification of (3.25) goes on a line similar to the one used for 
checking Eqn.(3.7) in the proof of Proposition 3.6. The left-hand side of
(3.25) is evaluated as:
\[
\frac{1}{2} \varphi ( \ 
(X)_{i_{n}i_{0}} (A)_{i_{0}i_{1}} \cdots (A)_{i_{n-1}i_{n}} \ )
\]
\[
= \ \frac{1}{2}  \delta_{i_{0}{\barl{i}}_{n}} \cdot
\delta_{i_{0}{\barl{i}}_{1}} \cdots \delta_{i_{n-1}{\barl{i}}_{n}} \cdot
\varphi ( \ (X)_{ {\barl{i}}_{0} i_{0}} 
(A)_{ i_{0} {\barl{i}}_{0} } \cdots  (A)_{ i_{n-1} {\barl{i}}_{n-1} }  \ )
\]
\begin{equation}
= \ \left\{  \begin{array}{lll}
2^{-1} \varphi( \xi a (a^{*}a)^{k} ),  & \mbox{if $n=2k+1$ and}  &
(i_{0},i_{1}, \ldots , i_{n} ) = (1,2, \ldots ,1,2)                   \\
                                       &                         &    \\
2^{-1} \varphi( \xi^{*} a^{*} (aa^{*})^{k} ), & \mbox{if $n=2k+1$ and}  &
(i_{0},i_{1}, \ldots , i_{n} ) = (2,1, \ldots ,2,1)                   \\
                                       &                         &    \\
0,                                     & \mbox{otherwise.}       &
\end{array} \right.
\end{equation}

The general term (indexed by $1 \leq m \leq n$) on the right-hand side 
of (3.25) is:
\[
\varphi_{2} ( \ \etaz ( \expd ( V_{ i_{0}i_{0} } A \cdots 
A V_{ i_{m-1}i_{m-1} } )) \cdot V_{ i_{m}i_{m} } A \cdots
A V_{ i_{n}i_{n} } \ )
\]
\[
= \ \varphi_{2} ( \ \etaz ( \delta_{ i_{0}i_{m-1} } \varphi ( 
(A)_{ i_{0}i_{1} }  \cdots  (A)_{ i_{m-2}i_{m-1} } ) V_{ i_{0}i_{0} } ) 
\cdot V_{ i_{m}i_{m} } A \cdots A V_{ i_{n}i_{n} } \ )
\]
\[
= \ \delta_{ i_{0}i_{m-1} } \varphi ( (A)_{ i_{0}i_{1} }  \cdots 
(A)_{ i_{m-2}i_{m-1} } ) \cdot \varphi_{2} ( \
V_{ {\barl{i}}_{0} {\barl{i}}_{0} }
V_{ i_{m}i_{m} } A \cdots A V_{ i_{n}i_{n} } \ )
\]
\[
= \ \frac{1}{2} \delta_{ i_{0}i_{m-1} } \cdot
\delta_{ i_{0}{\barl{i}}_{m} } \cdot
\delta_{ i_{0}{\barl{i}}_{1} } 
\delta_{ i_{1}{\barl{i}}_{2} } \cdots
\delta_{ i_{n-1}{\barl{i}}_{n} } \cdot
\varphi( (A)_{ i_{0}{\barl{i}}_{0} } \cdots
(A)_{ i_{m-2}{\barl{i}}_{m-2} } ) \cdot \varphi
( (A)_{ i_{m}{\barl{i}}_{m} } \cdots
(A)_{ i_{n-1}{\barl{i}}_{n-1} } ) 
\]
\begin{equation}
= \ \left\{  \begin{array}{ll}
2^{-1} \varphi( (aa^{*})^{(m-1)/2} ) \varphi( (a^{*}a)^{(n-m)/2} ), 
    & \mbox{if $m,n$ are both odd}                                     \\
    & \mbox{and } (i_{0},i_{1}, \ldots , i_{n} ) = (1,2, \ldots ,1,2)  \\ 
    &                                                                  \\
2^{-1} \varphi( (a^{*}a)^{(m-1)/2} ) \varphi( (aa^{*})^{(n-m)/2} ), 
    & \mbox{if $m,n$ are both odd}                                     \\
    & \mbox{and } (i_{0},i_{1}, \ldots , i_{n} ) = (2,1, \ldots ,2,1)  \\ 
    &                                                                  \\
0,  & \mbox{otherwise.}
\end{array} \right.
\end{equation}

By comparing (3.26) with (3.27) (and by also taking (3.24) into account)
we see that all it takes in order to obtain (3.25) is:
\begin{equation}
\varphi ( \xi a (a^{*}a)^{k} ) \ = \ 
\sum_{l=0}^{k} \ff((aa^*)^{l}) \cdot \ff ((a^*a)^{k-l}), \ \ 
\forall k \geq 0.
\end{equation}
Finally, we obtain (3.28) by using ``the odd half'' of the conjugate
relations (with respect to the scalars), which are fulfilled by $X$:
\begin{equation}
\ff(XA^{2k+1})  =  \sum_{l=1}^{2k+1} \ff_2(A^{l-1})\cdot
\ff_2 (A^{2k+1-l}) , \ \ k \geq 0.
\end{equation}
Indeed, we have:
\begin{eqnarray}
\ff_2(XA^{2k+1}) & = &
\ff_2
\left(  \begin{array}{cc}
\xi^*a^*(aa^*)^k  &   0    \\
0 & \xi a(a^*a)^k
\end{array}   \right)  \nonumber  \\
& = & \frac{1}{2} ( \ \varphi ( \xi^{*} a^{*} (aa^{*})^{k} ) +
     \ff(\xi a(a^*a)^k \ )   \nonumber \\
& = & \ff(\xi a(a^*a)^k), \mbox{ by (3.24);}  \nonumber
\end{eqnarray}
while on the other hand it is immediate that:
\begin{eqnarray}
\sum_{l=1}^{2k+1} \ff_2(A^{l-1})\cdot\ff_2 (A^{2k+1-l}) & = &
\sum_{l=0}^{k} \ff_2(A^{2l})\cdot\ff_2 (A^{2(k-l)})
\nonumber \\
& = & \sum_{l=0}^{k} \ff((aa^*)^{l})\cdot\ff ((a^*a)^{k-l}) \nonumber
\end{eqnarray}
(due to the particular form of $A$). So actually (3.28) reduces to 
(3.29). {\bf QED}

$\ $

We now discuss the special property of the $R$-diagonal element which 
will ensure the equality in (3.1) of Theorem 3.2.

$\ $

{\bf 3.8 Proposition.} In the framework of the Notations 3.4, we have
that: $Alg( \{ A \} \cup \D )$ is free from $\MM$ with amalgamation over
$\D$ if and only if $a$ is $R$-diagonal.

$\ $

In the proof of the Proposition 3.8 we will use the following lemma.

$\ $

{\bf 3.9 Lemma.} Let $( {\cal M} , \psi )$ be a $W^{*}$-probability 
space, and let us denote, for every $b \in {\cal M}$ and every $k \geq 1$:
\begin{equation}
\left\{ \begin{array}{lcl}
w_{11;k} (b)   &  = & (bb^{*})^{k} - \psi ( \ (bb^{*})^{k} \ ) I \\
w_{12;k} (b)   &  = & b(b^{*}b)^{k-1}                            \\
w_{21;k} (b)   &  = & b^{*} (bb^{*})^{k-1}                       \\
w_{22;k} (b)   &  = &  ( b^{*}b )^{k} - \psi ( \ ( b^{*}b )^{k} \ ) I.
\end{array}  \right.
\end{equation}
Then the following statements about an element $b \in {\cal M}$ are 
equivalent:

\vspace{4pt}

$1^{o}$ $b$ is $R$-diagonal in $( {\cal M} , \psi )$.

\vspace{4pt}

$2^{o}$ We have that:
\begin{equation}
\psi ( \  w_{ {\barl{i}}_{0} i_{1};k_{1} } (b)
w_{ {\barl{i}}_{1} i_{2};k_{2} } (b) \cdots
w_{ {\barl{i}}_{n-1} i_{n};k_{n} } (b)  \ ) \ = \ 0
\end{equation}
for every $n \geq 1$, $i_{0},i_{1}, \ldots , i_{n} \in \{ 1,2 \}$ and
$k_{1}, \ldots , k_{n} \geq 1$. (Same as in the preceding propositions,
we used in Equation (3.31) the convention of notation $\barl{i} = 3-i$,
for $i \in \{ 1,2 \}$.)

$\ $

{\bf Proof of Lemma 3.9.} $1^{o} \Rightarrow 2^{o}$ We can assume that 
$b =up$ where $u \in {\cal M}$ is a Haar unitary, $p=p^{*}$ is even
(i.e. $\psi ( p^{k} ) = 0$ for $k$ odd) and $\{ u,u^{*} \}$ is free 
from $\{ p \}$. Thus we have, for every $k \geq 1$:
\[
w_{11;k} (b) \ = \ u( p^{2k} - \psi ( p^{2k} ) I )u^{*} ,
\ \ w_{12;k} (b) \ = \ u p^{2k-1}, \ 
\]
\[
w_{21;k} (b) \ = \ p^{2k-1} u^{*}, \ \ 
w_{22;k} (b) \ = \ p^{2k} - \psi ( p^{2k} ) I.
\]
Every $w_{ij;k} (b)$ can be viewed as a word with 1, 2, or 3 letters
over the alphabet:
\begin{equation}
\{ u,u^{*} \} \cup \{ p^{k} - \psi (p^{k}) I \ | \ k \geq 1 \} ;
\end{equation}
and moreover the letters which form $w_{ij;k} (b)$ always come 
alternatively from $\{ u,u^{*} \}$ and
$\{ p^{k} - \psi (p^{k}) I \ | \ k \geq 1 \}$.

Given any $n \geq 1$, $i_{0},i_{1}, \ldots , i_{n} \in \{ 1,2 \}$ and
$k_{1}, \ldots , k_{n} \geq 1$, we claim that the product:
\begin{equation}
w \ := \ 
w_{ {\barl{i}}_{0} i_{1};k_{1} } (b)
w_{ {\barl{i}}_{1} i_{2};k_{2} } (b) \cdots
w_{ {\barl{i}}_{n-1} i_{n};k_{n} } (b)   
\end{equation}
still has the same alternance property of the letters, when viewed as a 
word over the alphabet (3.32). Indeed, for every $1 \leq m \leq n-1$ there
are two possibilities: either $i_{m} =1$, in which case
$w_{ {\barl{i}}_{m-1} i_{m};k_{m} } (b)$ ends with $u^{*}$ and 
$w_{ {\barl{i}}_{m} i_{m+1};k_{m+1} } (b)$ begins with a 
$p^{k} - \psi (p^{k}) I$; or $i_{m} = 2$, in which case
$w_{ {\barl{i}}_{m-1} i_{m};k_{m} } (b)$ ends with a 
$p^{k} - \psi (p^{k}) I$ and $w_{ {\barl{i}}_{m} i_{m+1};k_{m+1} } (b)$
begins with $u$. In both cases, the concatenation of
$w_{ {\barl{i}}_{m-1} i_{m};k_{m} } (b)$ and 
$w_{ {\barl{i}}_{m} i_{m+1};k_{m+1} } (b)$ is still alternating.

But if the product $w$ appearing in (3.33) is alternating when 
viewed as a word with letters from (3.32), then the equality 
$\psi (w) = 0$ follows from the definition of freeness (since
every letter in (3.32) is in the kernel of $\psi$, and since 
$\{ u,u^{*} \}$ is free from $\{ p \}$). 

\vspace{10pt}

$2^{o} \Rightarrow 1^{o}$. By enlarging the space $\ncpsm$ if necessary,
we can assume that there exists an $R$-diagonal element $c \in {\cal M}$,
such that $c^{*}c$ has the same distribution as $b^{*}b$. We denote 
$b_{1} := b$, $b_{2} := b^{*}$, $c_{1} := c$, $c_{2} := c^{*}$. We will 
show that:
\begin{equation}
\psi ( b_{i_{1}} b_{i_{2}} \cdots b_{i_{n}} ) \ = \ 
\psi ( c_{i_{1}} c_{i_{2}} \cdots c_{i_{n}} ), \ \ \forall n \geq 1, \
\forall i_{1}, \ldots , i_{n} \in \{ 1,2 \} . 
\end{equation}
From (3.34) it will follow that $b$ is $R$-diagonal (since $c$ is so, and
(3.34) means that $b$ and $c$ have the same $*$-distribution).

If $w_{ij;k} (c) \in {\cal M}$ is defined by analogy with Equation (3.30),
for $i,j \in \{ 1,2 \}$ and $k \geq 1$, then the implication
$1^{o} \Rightarrow 2^{o}$ proved above ensures that:
\begin{equation}
\psi ( \  w_{ {\barl{i}}_{0} i_{1};k_{1} } (c)
w_{ {\barl{i}}_{1} i_{2};k_{2} } (c) \cdots
w_{ {\barl{i}}_{n-1} i_{n};k_{n} } (c)  \ ) \ = \ 0
\end{equation}
for every $n \geq 1$, $i_{0},i_{1}, \ldots , i_{n} \in \{ 1,2 \}$ and
$k_{1}, \ldots , k_{n} \geq 1$. The equality in (3.34) will be obtained
by exploiting the similarity between (3.31) and (3.35).

We will prove (3.34) by induction on $n$. For $n=1$ we have to show that
$\psi (b) = \psi (c)$, $\psi ( b^{*} ) = \psi ( c^{*} )$. And indeed:
\[
\psi (b) \ = \ \psi ( w_{12;1} (b) ) \ \stackrel{(3.31)}{=} \ 0 
\ \stackrel{(3.35)}{=} \ \psi ( w_{12;1} (c) ) \ = \ \psi (c),
\]
while $\psi ( b^{*} ) = 0 = \psi ( c^{*} )$ can be shown in a similar way.

We consider now an $n \geq 2$. We assume that (3.34) is true for 
$1,2, \ldots , n-1$ and we prove it for $n$. Let us fix some indices 
$i_{1}, \ldots ,i_{n} \in \{ 1,2 \}$, about which we want to prove that 
(3.34) holds.

We take the product $b_{i_{1}} b_{i_{2}} \cdots b_{i_{n}}$, and draw a 
vertical bar between $b_{i_{m}}$ and $b_{i_{m+1}}$ for every 
$1 \leq m \leq n-1$ such that $i_{m} = i_{m+1}$. (For instance if
$b_{i_{1}} b_{i_{2}} \cdots b_{i_{n}}$ were to be $bb^{*}bbbb^{*}b^{*}b$,
then our bars would look like this: $bb^{*}b \mid b \mid bb^{*} \mid b^{*}b$.) 
By examining the sub-products of $b_{i_{1}} b_{i_{2}} \cdots b_{i_{n}}$ 
which sit between consecutive vertical bars, we find that we have written:
\begin{equation}
b_{i_{1}} b_{i_{2}} \cdots b_{i_{n}} \ = \ 
( w_{ {\barl{j}}_{0} j_{1};k_{1} } (b) + \lambda_{1} I ) \cdots 
( w_{ {\barl{j}}_{s-1} j_{s};k_{s} } (b) + \lambda_{s} I )
\end{equation}
for some $s \geq 1$, $j_{0},j_{1}, \ldots , j_{s} \in \{ 1,2 \}$,
$k_{1}, \ldots , k_{s} \geq 1$ having $k_{1} + \cdots + k_{s} = n$,
and $\lambda_{1}, \ldots , \lambda_{s} \in \C$. The number $\lambda_{r}$,
$1 \leq r \leq s$, is determined as follows: if $j_{r-1} = j_{r}$, then 
$\lambda_{r} = 0$; and if $j_{r-1} \neq j_{r}$, then 
$\lambda_{r} = \psi ( \ (b^{*}b)^{k_{r}} \ )$.

In a similar way we can write:
\begin{equation}
c_{i_{1}} c_{i_{2}} \cdots c_{i_{n}} \ = \ 
( w_{ {\barl{j}}_{0} j_{1};k_{1} } (c) + \lambda_{1} I ) \cdots 
( w_{ {\barl{j}}_{s-1} j_{s};k_{s} } (c) + \lambda_{s} I );
\end{equation}
and moreover, the parameters $s$, $j_{0},j_{1}, \ldots , j_{s}$,
$k_{1}, \ldots , k_{s}$, $\lambda_{1}, \ldots , \lambda_{s}$ appearing in 
(3.37) coincide with those from (3.36). Indeed, the values of $s$, 
$j_{0},j_{1}, \ldots , j_{s}$, $k_{1}, \ldots , k_{s}$ are determined
solely by how the vertical bars are placed between the $c_{i_{m}}$'s in 
$c_{i_{1}} c_{i_{2}} \cdots c_{i_{n}}$, and this is identical to how 
the vertical bars were placed in $b_{i_{1}} b_{i_{2}} \cdots b_{i_{n}}$.
After that, the value of every $\lambda_{r}$ is determined as 
$\delta_{ {\barl{j}}_{r-1}, j_{r} } \psi ( \ (c^{*}c)^{k_{r}} \ )$, which
is again the same as in (3.36), due to the fact that $b^{*}b$ and 
$c^{*}c$ have the same distribution.

By applying $\psi$ on both sides of (3.36) and then by expanding the 
product on the right-hand side, we obtain:
\[
\psi ( \ b_{i_{1}} b_{i_{2}} \cdots b_{i_{n}} \ ) \ = \ 
\psi ( \ w_{ {\barl{j}}_{0} j_{1};k_{1} } (b) \cdots 
w_{ {\barl{j}}_{s-1} j_{s};k_{s} } (b) \ )
\]
\[
+ \ \sum_{\emptyset \neq A \subseteq \{ 1, \ldots , s \} } \ 
\Bigl( \prod_{r \in A} \lambda_{r} \Bigr) \cdot \psi 
\Bigl( \prod_{r \in \{ 1, \ldots , s \} \setminus A}  
w_{ {\barl{j}}_{r-1} j_{r}; k_{r} } (b) \Bigr) .
\] 
The corresponding operations done in (3.37) yield an identical formula,
where we have $c$'s instead of $b$'s. But we know that 
\[
\psi ( \  w_{ {\barl{j}}_{0} j_{1};k_{1} } (b) \cdots
w_{ {\barl{j}}_{s-1} j_{s};k_{s} } (b)  \ )
\ \stackrel{(3.31)}{=} \ 0 \ \stackrel{(3.35)}{=} \ 
\psi ( \  w_{ {\barl{j}}_{0} j_{1};k_{1} } (c) \cdots
w_{ {\barl{j}}_{s-1} j_{s};k_{s} } (c)  \ ),
\]
while on the other hand the induction hypothesis gives us that:
\[
\psi \Bigl( \prod_{r \in \{ 1, \ldots , s \} \setminus A}  
w_{ {\barl{j}}_{r-1} j_{r}; k_{r} } (b) \Bigr) \ = \ 
\psi \Bigl( \prod_{r \in \{ 1, \ldots , s \} \setminus A}  
w_{ {\barl{j}}_{r-1} j_{r}; k_{r} } (c) \Bigr) ,
\]
for every $\emptyset \neq A \subseteq \{ 1, \ldots ,s \}$. These 
equalities imply in turn that
$\psi ( \ b_{i_{1}} b_{i_{2}} \cdots b_{i_{n}} \ )$ = 
$\psi ( \ c_{i_{1}} c_{i_{2}} \cdots c_{i_{n}} \ )$, as desired.
{\bf QED}

$\ $

{\bf Proof of Proposition 3.8.} For every $i,j \in \{ 1,2 \}$ and
$k \geq 1$ we denote by $w_{ij;k} (a)$ the element of $\A$ defined by the 
same recipe as in Equation (3.30) of Lemma 3.9, and we denote by
$W_{ij;k} \in M_{2} ( \A )$ the matrix which has its $(i,j)$-entry equal 
to $w_{ij;k} (a)$, and its other entries equal to 0. 

It is immediately seen that $Alg( \{ A \} \cup \D )$ is linearly spanned 
by the matrices of the form:
\[
\left(  \begin{array}{cc} 
(aa^{*})^{k}  &  0   \\  0   &  0
\end{array}  \right),  \ \
\left(  \begin{array}{cc} 
0  & a(a^{*}a)^{k}   \\  0   &  0
\end{array}  \right),  \ \
\left(  \begin{array}{cc} 
0  &  0   \\  a^{*}(aa^{*})^{k}  &  0   
\end{array}  \right),  \ \
\left(  \begin{array}{cc} 
0  &  0  \\  0  & (a^{*}a)^{k}   
\end{array}  \right),  \ \  k \geq 0;
\]
and this implies the formula:
\begin{equation}
\{ X \in Alg( \{ A \} \cup \D )  \ | \ \expd (X) = 0 \} \ = \
\mbox{span} \{ W_{ij;k} \ | \ i,j \in \{ 1,2 \} , \ k \geq 1 \} .
\end{equation}
On the other hand it is clear that:
\begin{equation}
\{ X \in \MM \ | \ \expd (X) = 0 \} \ = \
\mbox{span} \{ V_{12} , V_{21}  \} .
\end{equation}

From (3.38-39) it follows that $Alg( \{ A \} \cup \D )$ is free from
$\MM$ with amalgamation over $\D$ if and only if:
\begin{equation}
\left\{ \begin{array}{l}
\expd ( \ U' W_{ j_{1}' j_{1}'';k_{1} } V_{ i_{1} {\barl{i}}_{1} }  \cdots 
V_{ i_{n-1} {\barl{i}}_{n-1} } W_{ j_{n}' j_{n}''; k_{n} } U'' \ ) 
\ = \ 0,                  \\
                          \\
\forall n \geq 1, \
\forall j_{1}' , j_{1}'', \ldots , j_{n}', j_{n}'', i_{1}, \ldots ,
i_{n-1} \in \{ 1,2 \} ,      \\
\forall k_{1} , \ldots , k_{n} \geq 1, \ \forall U', U'' \in 
\{ V_{11} + V_{22}, V_{12}, V_{21} \} .
\end{array}  \right.
\end{equation}
The matrix product appearing in (3.40) is 0 if it is not true that 
$j_{1}'' = i_{1}$, ${\barl{i}}_{1} = j_{2}' , \ldots , j_{n-1}'' = i_{n-1}$,
${\barl{i}}_{n-1} = j_{n}'$. And consequently, (3.40) is equivalent to:
\begin{equation}
\left\{ \begin{array}{l}
\expd ( \ U' W_{ {\barl{i}}_{0} i_{1};k_{1} } V_{ i_{1} {\barl{i}}_{1} } 
\cdots W_{ {\barl{i}}_{n-2} i_{n-1};k_{n-1} } 
V_{ i_{n-1} {\barl{i}}_{n-1} }
W_{ {\barl{i}}_{n-1} i_{n};k_{n} } U'' \ ) \ = \ 0,                  \\
                                                                   \\
\forall n \geq 1, \
\forall i_{0}, i_{1}, \ldots , i_{n} \in \{ 1,2 \} ,               \\
\forall k_{1} , \ldots , k_{n} \geq 1, \ \forall U', U'' \in 
\{ V_{11} + V_{22}, V_{12}, V_{21} \} .
\end{array}  \right.
\end{equation}
But now, the matrix product appearing in (3.41) has one entry 
equal to $w_{ {\barl{i}}_{0} i_{1};k_{1} } (a)$ $\cdots$
$w_{ {\barl{i}}_{n-1} i_{n};k_{n} } (a)$ (which can appear on any of 
the four possible positions, depending on the choices of $U'$ and $U''$); 
and has the other three entries equal to 0. This makes it immediate that 
the condition (3.41) is equivalent to the one presented in the 
Lemma $3.9.2^{o}$ (applied here to $\ncps$). {\bf QED}

$\ $

{\bf Proof of Theorem 3.2.} The inequality (3.1) is obtained by 
putting together the Equations (3.5) and (3.18) established in the 
Propositions 3.6, 3.7. From (3.5) and (3.18) it is also clear that
(3.1) holds with equality if and only if:
\begin{equation}
\Phi^{*} (A: \MM , \eta ) \ = \ \Phi^{*} (A: \D , \etaz ).
\end{equation}
As reviewed in the Remark $2.8.3^{o}$, a sufficient condition for (3.42)
to take place is that $Alg( \{ A \} \cup \D )$ and $\MM$ are free with 
amalgamation over $\D$. But by Proposition 3.8, this sufficient 
condition is equivalent to the fact that $a$ is an $R$-diagonal element.
{\bf QED}

$\ $

{\bf 3.10 Remark.} If one is only interested in establishing the 
inequality (3.1), then a substantial short-cut can be taken through
the above considerations. The short-cut goes by verifying directly
that if $( \xi , \xi^{*} )$ is a conjugate system for $( a,a^{*} )$, 
then $X := \left( \begin{array}{cc} 0  & \xi^{*} \\ \xi & 0
\end{array} \right)$ fulfills the conjugate relations for 
$A := \left( \begin{array}{cc} 0  & a \\ a^{*} & 0 \end{array} \right)$, 
with respect to the scalars; this immediately implies 
$\Phi^{*} (a,a^{*} ) \geq 2 \Phi^{*} (A) = 2 \Phi^{*} ( \mu )$, 
i.e. (3.1).

The reason for insisting to put into evidence the relations shown in 
(3.5) and (3.18) is that they also give us a non-trivial necessary and
sufficient condition -- Equation (3.42) -- for the minimal free
Fisher information to be attained. As mentioned in Remark 3.3, the problem
of determining what are the $*$-distributions of $a$ which attain the 
minimal $\Phi^{*} (a,a^{*})$ is open, and seemingly difficult. The 
condition in (3.42) helps clarifying the nature of this problem, by
reducing it to the following one (also open):

$\ $

{\bf 3.11 Problem.} Is it true that if 
$\Phi^{*} (A: \MM , \eta ) = \Phi^{*} ( A: \D , \etaz ) < \infty$, 
then necessarily $Alg( \{ A \} \cup \D )$ and $\MM$ are free with 
amalgamation over $\D$?

$\ $

Settling the Problem 3.11 in the affirmative would imply that if the 
minimal value of $\Phi^{*} ( a,a^{*} )$ (under the constraint in (1.1)) 
is finite, then this minimal value can be reached only by an $R$-diagonal
element. While on the other hand, a negative answer in 3.11 would provide 
examples of situations when the infimum in (1.1) is finite and is reached
by non-$R$-diagonal elements.

$\ $

$\ $

$\ $

\setcounter{section}{4}
\setcounter{equation}{0}
{\large\bf 4. Minimization of free Fisher information for matrix entries}

$\ $

Let $d$ be a fixed positive integer. We will consider here the following
two minimization problems:

\vspace{6pt}

(a) Determine the minimal possible value of $\Phi^{*} ( \ \aijaijs \ )$, 
if the family 
\newline
$\aijs$ (of elements in some $W^{*}$-probability space
$\ncps$, with $\varphi$ faithful trace) is such that the matrix 
$A = \mataij$ has a prescribed $*$-distribution.

\vspace{6pt}

(b) Determine the minimal possible value of $\Phi^{*} ( \ \bijs \ )$, 
if the family $\bijs$
(of elements in some $W^{*}$-probability space $\ncps$, with $\varphi$
faithful trace) is such that the matrix $B = \matbij$ is selfadjoint, 
and has a prescribed distribution. Note that if $B=B^{*}$, then $\bijs$ is 
a selfadjoint family of elements of $\A$ -- indeed, the involution
$\sigma (i,j) := (j,i)$ has the property that
$b_{i,j}^{*} = b_{\sigma (i,j)}$, for every $1 \leq i,j \leq d$.

$\ $

The solutions of these problems are provided by the Theorem 1.2 stated in
the Introduction. For instance for (b) we have that, given a probability 
measure $\mu$ with compact support on $\R$:
\begin{equation}
\min \Bigl\{ \ \Phi^* ( \ \bijs \ ) \ \mid \ 
\begin{array}{c}
B = \matbij = B^{*}   \\
\mbox{ has distribution } \mu
\end{array}  \Bigr\} \ = \ d^{3} \Phi^{*} ( \mu ).
\end{equation}
A similar formula holds in the framework of the problem (a) (but where 
for the role of $\mu$ we must now consider a linear functional on
$\C \langle X,X^{*} \rangle$, which can appear as a $*$-distribution 
in a tracial $W^{*}$-probability space).

In order to infer (4.1) as a consequence of Theorem $1.2.2^{o}$ (and the
corresponding conclusion from Theorem $1.2.1^{o}$), there is one more
detail that needs to be verified -- that the freeness conditions appearing
in Theorem 1.2 can indeed be fulfilled, in the context where the joint
distribution of $A$ and $A^{*}$ (in $1^{o}$) and the distribution of $B$
(in $2^{o}$) are prescribed. We discuss in more detail the selfadjoint 
case of $2^{o}$; the non-selfadjoint case is similar.

So, let $\mu$ be a fixed probability measure with compact
support on $\R$. One can find a $W^{*}$-probability space $\ncpsm$, with 
$\psi$ faithful trace, and $x$ and $\vijs$ in $\M$ such that:
\begin{itemize}
\item[(i)]
$x = x^{*}$ has distribution $\mu$;
\item[(ii)]
the $v_{ij}$'s form a family of matrix units (i.e., $v_{ij}v_{kl} =
\delta_{jk}v_{il}$, $v_{ij}^{*} = v_{ji}$, $\forall 1 \leq i,j,k,l \leq d$,
and $\sum_{i=1}^{d} v_{ii} =I$);
\item[(iii)]
$x$ is free from $\vijs$.
\end{itemize}
(An example of such $\ncpsm$ is the free product $( L^{\infty} ( \mu ) , 
d \mu ) \star ( M_{d} ( {\bf C} ), tr)$. ) Consider the 
compressed $W^{*}$-probability space $\ncps$, where
$\A := v_{11} \M v_{11}$ and $\varphi ( \cdot ) := d \psi ( \cdot )$ 
on $\A$; and in $\A$ consider the family of compressions
$b_{ij} \ := \ v_{1i} x v_{j1}, \ \ 1 \leq i,j \leq d$.
Then the self-adjoint matrix $B = \matbij$ has distribution $\mu$  in 
$\ncpsd$, and is on the other hand free from $M_{d} ( \C I) \subseteq
M_{d} ( \A )$. These things happen because 
the spaces $\ncpsm$ and $\ncpsd$ are isomorphic, via the $*$-isomorphism
$\M \ni y \mapsto ( v_{1i}yv_{j1} )_{i,j=1}^{d} \in  M_{d} ( \A )$,
which sends $x$ to $B$ and span$\{ v_{ij} \ | \ 1 \leq i,j \leq d \}$ 
onto $M_{d} ( \C I )$.

$\ $

It thus remains that we prove the Theorem 1.2. We will in fact prove 
more, namely:

$\ $

{\bf 4.1 Proposition.} Let $\ncps$ be a $W^{*}$-probability space, with
$\varphi$ faithful trace, and let $\B \subseteq \A$ be a unital 
$W^{*}$-subalgebra. Let $d$ be a positive integer; consider the 
$W^{*}$-probability space $\ncpsd$ (defined as in Notations $2.1.2^{o}$),
and the $W^{*}$-subalgebra $M_{d} ( \B ) \subseteq M_{d} ( \A )$.

\vspace{4pt}

$1^{o}$ For every $A = \mataij \in M_{d} ( \A )$, we have:
\begin{equation}
\Phi^{*} ( \ \aijaijs : \B \ ) \ = \
d^{3} \Phi^{*} ( \ \{ A, A^{*} \} : M_{d} ( \B ) \ ).
\end{equation}

\vspace{4pt}

$2^{o}$ For every $G = ( g_{ij} )_{i,j=1}^{d} \in M_{d} ( \A )$ such that
$G = G^{*}$, we have:
\begin{equation}
\Phi^{*} ( \ \{ g_{ij} \}_{1 \leq i,j \leq d} : \B \ ) \ = \
d^{3} \Phi^{*} ( \ G: M_{d} ( \B ) \ ).
\end{equation}

$\ $

The name of the selfadjoint matrix appearing in $4.1.2^{o}$ was changed 
to $G$ (from $B$, as was in Theorem 1.2) in order to avoid any confusion
with the elements of $M_{d} ( \B )$. Note that if in Proposition 4.1
we take $\B = \C I$, then the Equations (4.2) and (4.3) become:
\begin{equation}
\Phi^{*} ( \ \aijaijs \ ) \ = \
d^{3} \Phi^{*} ( \ \{ A, A^{*} \} : M_{d} ( \C I ) \ ),
\end{equation}
and respectively
\begin{equation}
\Phi^{*} ( \ \{ g_{ij} \}_{1 \leq i,j \leq d} \ ) \ = \
d^{3} \Phi^{*} ( \ G: M_{d} ( \C I ) \ ).
\end{equation}
The statements of Theorem 1.2 follow immediately from these relations.
Indeed, for $1.2.1^{o}$ we only have to use (4.4) and the fact (reviewed
in Remark 2.6) that $\Phi^{*} ( \ \{ A , A^{*} \} : M_{d} ( \C I ) \ )
\geq \Phi^{*} (A,A^{*})$, with equality when $\{ A,A^{*} \}$ is free 
from $M_{d} ( \C I )$; and similarly for $1.2.2^{o}$.

$\ $

{\bf Proof of Proposition 4.1.} The proofs of $4.1.1^{o}$ and $4.1.2^{o}$
are similar to each other (and also similar to the proof of Proposition 
3.6 from the previous section). For this reason, we will only do 
$4.1.1^{o}$, and  leave $4.1.2^{o}$ as an exercise to the reader.

In $4.1.1^{o}$ we first consider the situation when
$\Phi^{*} ( \ \aijaijs : \B \ ) < \infty$. In this case, the family
$\aijaijs$ has a conjugate system $\xiijxiijs$ with respect to $\B$.
Let us define:
\begin{equation}
X \ := \ \frac{1}{d} ( \xi_{ji} )_{i,j=1}^{d} \in \ldaphi
\end{equation}
(where the identification discussed in $2.1.3^{o}$ is used). We will
show that $\{ X, X^{*} \}$ is a conjugate system for $\{ A,A^{*} \}$,
with respect to $M_{d} ( \B )$. This will entail (4.2) (under the 
hypothesis $\Phi^{*} ( \ \aijaijs \ ) < \infty$), because it will give:
\[
\Phi^{*} ( \ \{ A,A^{*} \} \ : M_{d} ( \B ) \ ) \ = \
|| X ||_{\lphid}^{2} + || X^{*} ||_{\lphid}^{2}
\]
\[
\stackrel{(2.3)}{=} \ \frac{1}{d} \sum_{i,j=1}^{d}
( \ || \frac{1}{d} \xi_{ji} ||_{\lphi}^{2} +
|| \frac{1}{d} \xi_{ij}^{*} ||_{\lphi}^{2}  \ )
\]
\[
= \ \frac{1}{d^3} \sum_{i,j=1}^{d} ( \
|| \xi_{ij} ||_{\lphi}^{2} + || \xi_{ij}^{*} ||_{\lphi}^{2} \ )
\ = \ \frac{1}{d^3} \Phi^{*} ( \ \aijaijs \ ).
\]

The conjugate relations which we need to verify are:
\begin{equation}
\varphi_{d} ( XB_{0}A_{i_{1}}B _{1} \cdots A_{i_{n}} B_{n} ) \ = \
\end{equation}
\[
\sum_{m=1}^{n} \delta_{i_{m},1} \cdot  \varphi_{d} 
( B_{0}A_{i_{1}} \cdots A_{i_{m-1}} B_{m-1} ) \cdot 
\varphi_{d} ( B_{m}A_{i_{m+1}} \cdots A_{i_{n}} B_{n} ), 
\]
for $n \geq 1$, $B_{0},B_{1}, \ldots , B_{n} \in M_{d} ( \B )$ and 
$i_{1},i_{2}, \ldots , i_{n} \in \{ 1,2 \}$, where we denoted 
$A_{1} := A, A_{2} := A^{*}$. The list of conjugate relations for 
$\{ A,A^{*} \}$ with respect to $M_{d} ( \B )$ also contains:
\begin{equation}
\varphi_{d} ( X^{*} B_{0}A_{i_{1}}B _{1} \cdots A_{i_{n}} B_{n} ) \ = 
\end{equation}
\[
\sum_{m=1}^{n} \delta_{i_{m},2} \cdot  \varphi_{d} 
( B_{0}A_{i_{1}} \cdots A_{i_{m-1}} B_{m-1} ) \cdot 
\varphi_{d} ( B_{m}A_{i_{m+1}} \cdots A_{i_{n}} B_{n} ), 
\]
(for $n \geq 1$, $B_{0},B_{1}, \ldots , B_{n} \in M_{d} ( \B )$, 
$i_{1},i_{2}, \ldots , i_{n} \in \{ 1,2 \}$), and 
\begin{equation}
\varphi_{d} ( XB ) \ = \
\varphi_{d} ( X^{*}B ) \ = \ 0, \ \ \forall B \in M_{d} ( \B ) .
\end{equation}
But however, (4.8) readily follows from (4.7) by taking an adjoint 
and then doing a circular permutation under $\varphi_{d}$; while
(4.9) is a direct consequence of the equations 
$\varphi ( \xi_{ij} b) = 0$, $1 \leq i,j \leq d$, $b \in \B$, which 
appear on the list of conjugate relations satisfied by the family
$\xiijxiijs$. (Hence indeed, only (4.7) needs to be checked.)

For every $1 \leq k,l \leq d$ and every $b \in \B$ let us denote by
$V_{kl} \otimes b$ the matrix in $M_{d} ( \B )$ which has its 
$(k,l)$-entry equal to $b$, and all its other entries equal to 0.
By multilinearity we can assume in (4.7) that
$B_{0} = V_{k_{0}l_{0}} \otimes b_{0}, \ldots ,
B_{n} = V_{k_{n}l_{n}} \otimes b_{n}$ for some
$k_{0},l_{0} , \ldots , k_{n},l_{n} \in \{ 1, \ldots , d \}$
and $b_{0}, \ldots , b_{n} \in \B$. The left-hand side of (4.7) is
then equal to:
\[
\varphi_{d} ( \ (V_{l_{n}l_{n}} \otimes I)X 
(V_{ k_{0}l_{0} } \otimes b_{0}) A_{i_{1}} 
(V_{ k_{1}l_{1} } \otimes b_{1})  \cdots A_{i_{n}}
(V_{ k_{n}l_{n} } \otimes b_{n}) \ )
\]
\[
= \ \frac{1}{d} \varphi ( \ (X)_{l_{n}k_{0}} b_{0}
(A_{i_{1}})_{ l_{0}k_{1} } b_{1} \cdots
(A_{i_{n}})_{ l_{n-1}k_{n} } b_{n} \ )
\]
\begin{equation}
= \ \frac{1}{d^{2}} \varphi ( \ \xi_{k_{0}l_{n}} b_{0} 
(A_{i_{1}})_{ l_{0}k_{1} } b_{1} \cdots
(A_{i_{n}})_{ l_{n-1}k_{n} } b_{n} \ ),
\end{equation}
where ``$(A_{i_{1}})_{ l_{0}k_{1} }$'' stands for the 
$(l_{0},k_{1})$-entry of the matrix $A_{i_{1}}$, etc (same conventions
of notation as in Section 3). By using the conjugate relations 
satisfied by the family $\xiijxiijs$, we can continue (4.10) with:
\begin{equation} 
= \ \frac{1}{d^2} \sum_{m=1}^{n} \delta_{i_{m},1}
\delta_{k_{0},l_{m-1}} \delta_{l_{n},k_{m}} \cdot
\varphi ( \ b_{0} (A_{i_{1}})_{ l_{0}k_{1} } \cdots
(A_{i_{m-1}})_{ l_{m-2}k_{m-1} } b_{m-1} \ ) \cdot
\end{equation} 
\[
\cdot \varphi ( \ b_{m} (A_{i_{m+1}})_{ l_{m}k_{m+1} } \cdots
(A_{i_{n}})_{ l_{n-1}k_{n} } b_{n} \ ).
\]
It is straightforward to observe that the summation which appeared in
(4.11) is equal to the right-hand side of (4.7).

In order to complete the verification that $\{ X,X^{*} \}$ is the 
conjugate of $\{ A,A^{*} \}$ with respect to $M_{d} ( \B )$, we must 
also show that:
\begin{equation}
X \ \in {\overline{Alg( \{ A,A^{*} \} \cup M_{d} ( \B ) 
)}}^{|| \cdot ||_{2}}.
\end{equation}
For every $1 \leq i,j \leq d$ let us denote by $A_{ij} \in M_{d} ( \A )$
and respectively by $X_{ij} \in \ldaphi$ the matrix which has $a_{ij}$
(respectively $\xi_{ij}$) on its (1,1)-entry, and 0's on all the other
entries. Then:
\[
A_{ij} \ = \ (V_{1i} \otimes I) A (V_{j1} \otimes I)
\in Alg( \{ A,A^{*} \} \cup M_{d} ( \B ) ), \ \
\forall 1 \leq i,j \leq d.
\]
Among the properties satisfied by $\xiijxiijs$ (as conjugate for
$\{ a_{ij} , a_{ij}^{*} \}_{1 \leq i,j \leq d}$), we also have that:
\[
\xi_{kl} \in 
{\overline{Alg( \aijaijs  \cup \B )}}^{|| \cdot ||_{2}}, \ \
\forall 1 \leq k,l \leq d.
\]
Consequently, by forming polynomials with matrices of the form 
$A_{ij}, A_{ij}^{*}$ and $V_{11} \otimes b$ ($b \in \B$), and then
by taking $|| \cdot ||_{2}$-limits, we obtain that every 
$X_{kl}$ ($1 \leq k,l \leq d$) belongs to the $|| \cdot ||_{2}$-closed 
space indicated in (4.12). This space is invariant under the left/right 
action of elements from $M_{d} ( \B )$, hence we can conclude 
that it also contains
\[
X \ = \ \frac{1}{d} \sum_{k,l=1}^{d} ( V_{l1} \otimes I) X_{kl}
( V_{1k} \otimes I ), 
\]
as desired.

Equation (4.2) is now proved in the case when
$\Phi^{*} ( \ \aijaijs : \B \ ) < \infty$. It remains to show that
$\Phi^{*} ( \ \aijaijs : \B \ ) = \infty$ $\Rightarrow$
$\Phi^{*} ( \{ A,A^{*} \} \ : M_{d} ( \B ) ) = \infty$; or equivalently,
that $\Phi^{*} ( \{ A,A^{*} \} \ : M_{d} ( \B ) ) < \infty$ $\Rightarrow$
$\Phi^{*} ( \ \aijaijs : \B \ ) < \infty$.

If $\Phi^{*} ( \{ A,A^{*} \} \ : M_{d} ( \B ) ) < \infty$, then there 
exists $X \in \ldaphi$ such that $\{ X,X^{*} \}$ fulfills the conjugate 
relations for $\{ A,A^{*} \}$, with respect to $M_{d} ( \B )$. We write 
$X$ as a $d \times d$-matrix (as in $2.1.3^{o}$):
\[
X \ = \ ( \eta_{ij} )_{i,j=1}^{d}, \ \
\mbox{with } \eta_{ij} \in \laphi , \ 1 \leq i,j \leq d;
\]
and we set $\xi_{ij} :=d \eta_{ji}$, $1 \leq i,j \leq d$. We claim that
$\xiijxiijs$ fulfills the conjugate relations for $\aijaijs$, with respect
to $\B$. Since the calculation verifying this claim is very similar in
spirit with the one which concluded the proof of Proposition 3.6, we will 
only mention its guiding line, and leave the details to the reader. The 
generic relation that needs to be proved is of the form:
\begin{equation}
\varphi ( \ \xi_{kl} b_{0} (A_{i_{1}})_{k_{1}l_{1}} b_{1} \cdots
(A_{i_{n}})_{k_{n}l_{n}} b_{n}  \ ) \ = \ 
\sum_{m=1}^{n} \delta_{i_{m},1}
\delta_{k,k_{m}} \delta_{l,l_{m}} \cdot
\end{equation}
\[
\varphi ( \ b_{0} (A_{i_{1}})_{ k_{1}l_{1} } \cdots
(A_{i_{m-1}})_{ k_{m-1}l_{m-1} } b_{m-1} \ ) \cdot
\varphi ( \ b_{m} (A_{i_{m}})_{ k_{m}l_{m} } \cdots
(A_{i_{n}})_{ k_{n}l_{n} } b_{n} \ ),
\]
for $n \geq 1$, $b_{0}, \ldots , b_{n} \in \B$,
$i_{1}, \ldots , i_{n} \in \{ 1,2 \}$,
$k_{1}, l_{1}, \ldots , k_{n},l_{n} \in \{ 1, \ldots ,d \}$. The line 
for establishing (4.13) goes by writing its left-hand side as
\begin{equation}
d^{2} \varphi_{d} ( \ X ( V_{k,k_{1}} \otimes b_{0}) A_{i_{1}}
(V_{l_{1},k_{2}} \otimes b_{1} ) A_{i_{2}}
(V_{l_{2},k_{3}} \otimes b_{2} ) \cdots A_{i_{n}} 
(V_{l_{n},l} \otimes b_{n} )  \ );
\end{equation}
then by using in (4.14) the conjugate relations fulfilled by 
$\{ X,X^{*} \}$; and finally by evaluating (in a straightforward way)
the terms of the summation which is obtained in this manner.

But if $\xiijxiijs$ fulfills the conjugate relations for $\aijaijs$,
with respect to $\B$, then it follows that
$\Phi^{*} ( \ \aijaijs : \B \ ) < \infty$, and this concludes the proof.
\newline
{\bf QED}

$\ $

By using the Theorem 1.2, we can now easily prove the generalization of
our minimization result for $\Phi^{*}$, which was stated in Theorem 1.3.

$\ $

{\bf Proof of Theorem 1.3.} Let us fix a probability measure $\nu$ with 
compact support on $[0, \infty )$, and a positive integer $d$. We denote
the symmetric square root of $\nu$ (defined as in 3.1) by $\mu$.

Let $\ncps$ be a $W^{*}$-probability space, with $\varphi$ faithful 
trace, and let $\aijs$ be elements of $\A$ such that if we set 
$A:= \mataij$, then $A^{*}A$ has distribution $\nu$ in $\ncpsd$. Then:
\begin{equation}
\Phi^* ( A , A^* ) \ \geq \ 2 \Phi^* ( \mu ) 
\end{equation}
(by Theorem 1.1); if we combine this with the inequality (1.4) of 
Theorem 1.2, we get:
\begin{equation}
\Phi^{*} ( \ \aijaijs \ ) \ \geq \ 2 d^{3} \Phi^{*} ( \mu ).
\end{equation}
A discussion similar to the one preceding Proposition 4.1 shows that
(in the context where $\nu$ and $d$ are prescribed) we can pick the 
family $\aijs$ such that in addition to the condition that the 
distribution of $A^{*}A$ be $\nu$, we also have:
\begin{itemize}
\item[(i)]
$A$ is $R$-diagonal in $M_{d} ( \A )$; and
\item[(ii)]
$\{ A, A^{*} \}$ is free from the algebra of scalar matrices 
$M_{d} ( \C I ) \subseteq M_{d} ( \A )$.
\end{itemize}

The condition (i) implies that (4.15) holds with equality, while (ii) 
implies equality in (1.4) of Theorem 1.2; hence (i)+(ii) ensure that 
the lower bound $2 d^{3} \Phi^{*} ( \mu )$ of (4.16) is actually
attained.  {\bf QED}

$\ $

In the case when $\Phi^{*} ( \mu ) < \infty$, it would be interesting to 
know if the conditions (i) and (ii) mentioned in the proof of Theorem 1.3
are also necessary for (4.16) to hold with equality. Deciding on this
fact would amount to solving the Problem 3.11 (which corresponds to the 
particular case $d=1$), and another problem of the same nature -- whether
the equality $\Phi^{*} ( \{ A,A^{*} \} : M_{d} ( \C I ) )$ =
$\Phi^{*} ( A,A^{*} ) < \infty$ must imply the freeness of 
$\{ A,A^{*} \}$ and $M_{d} ( \C I )$.
 
$\ $

$\ $

$\ $

\setcounter{section}{5}
\setcounter{equation}{0}
{\large\bf 5. The corresponding maximization problems for the 
free entropy $\chi^{*}$ } 

$\ $

In this section we will consider the concept of free entropy $\chi^*$,
defined in \cite{V-P5} in terms of the free information $\Phi^*$.
We will treat the questions of maximizing $\chi^{*}$, under constraints
similar to those discussed in the previous sections. The results
concerning $\chi^*$ will follow from the corresponding results for the 
free Fisher information. 

Let $\ncps$ be a $W^{*}$-probability space with $\varphi$ a faithful
trace, and consider a selfadjoint family of elements of $\A$ which is 
given in the form: $\{ a_{i}, a_{i}^{*} \}_{1 \leq i \leq m} \cup
\{ b_{j} \}_{1 \leq j \leq n}$, where $b_{j} = b_{j}^{*}$ for $1 \leq j
\leq n$. By enlarging $\ncps$ if necessary, we can assume there exist 
circular elements $c_{1}, \ldots , c_{m} \in \A$ and 
semicircular elements $s_{1}, \ldots , s_{n} \in \A$ such that
$\{ c_{1}, c_{1}^{*} \} , \ldots , \{ c_{m}, c_{m}^{*} \}$, 
$\{ s_{1} \} , \ldots, \{ s_{n} \} ,
\{ a_{1}, a_{1}^{*} , \ldots , a_{m}, a_{m}^{*}, 
b_{1} , \ldots, b_{n} \}$ are free. We will assume in addition that
$c_{1}, \ldots , c_{m}$ and $s_{1}, \ldots , s_{n}$ are normalized by
their variance (i.e. $\varphi ( c_{i}^{*} c_{i} )$ = 1 =
$\varphi ( s_{j}^{2} )$, for every $1 \leq i \leq m$, $1 \leq j \leq n$).
Then the free entropy
$\chi^{*} ( \ \{ a_{i}, a_{i}^{*} \}_{1 \leq i \leq m} \cup
\{ b_{j} \}_{1 \leq j \leq n} \ )$ $\in$ $[ - \infty , \infty )$ is defined
by the formula
\begin{equation}
\chi^{*} ( \ \{ a_{i}, a_{i}^{*} \}_{1 \leq i \leq m} \cup
\{ b_{j} \}_{1 \leq j \leq n} \ ) \ = \ 
\frac{2m+n}{2} \log ( 2 \pi e ) +
\end{equation}
\[
+ \frac{1}{2} \int_{0}^{\infty} \Bigl( \ \frac{2m+n}{1+t} - 
\Phi^{*} ( \ \{ a_{i}+ \sqrt{t} c_{i}, a_{i}^{*} + \sqrt{t} c_{i}^{*}
\}_{1 \leq i \leq m} \cup
\{ b_{j} + \sqrt{t} s_{j} \}_{1 \leq j \leq n} \ ) \ \Bigr) dt .
\]
The integral on the right-hand side of (5.1) makes sense, and takes indeed
value in $[ - \infty , \infty )$ -- see Corollary 6.14, Proposition 7.2 
in \cite{V-P5}. (In order to apply literally the estimates from \cite{V-P5},
one first replaces every pair $\{ c_{j} , c_{j}^{*} \}$ with the pair of
selfadjoints $\{ (c_{j} + c_{j}^{*})/ \sqrt{2}$, 
$(c_{j} - c_{j}^{*})/ i \sqrt{2} \}$ -- this does not affect the integrand
on the right-hand side of (5.1).) Moreover, the value of the integral in
(5.1) does not depend on the choice of $c_{1}, \ldots , c_{m}$,
$s_{1}, \ldots , s_{n}$; in fact it is easy to see that 
$\chi^{*} ( \ \{ a_{i}, a_{i}^{*} \}_{1 \leq i \leq m} \cup
\{ b_{j} \}_{1 \leq j \leq n} \ )$ depends only on the joint distribution
of $\{ a_{1}, a_{1}^{*}, \ldots ,a_{m}, a_{m}^{*}$, 
$b_{1}, \ldots , b_{n} \}$ in $\ncps$.
 
If $\mu$ is a probability measure with compact support on $\R$, then we will
denote (similarly to how we did with $\Phi^{*}$ in Notation 2.9):
\begin{equation}
\chi^{*} ( \mu ) \ := \ \chi^{*} (x),
\end{equation}
where $x$ is an arbitrary selfadjoint random variable with distribution 
$\mu$. Similarly to the situation for $\Phi^{*}$, there exists an explicit 
integral formula for $\chi^{*} ( \mu )$, namely:
\begin{equation}
\chi^{*} ( \mu ) \ = \ \int \int  \log | \ s-t \ | \ d \mu (s)
d \mu (t) \ + \ \frac{3}{4} \ + \ \frac{\log (2 \pi )}{2}
\end{equation}
(\cite{V-P2}, Proposition 4.5, combined with \cite{V-P5}, Proposition 7.6).

$\ $

We now start towards the proofs of Theorems 1.4 and 1.5. Following 
the same line which we used for $\Phi^{*}$, we will first do the
Theorem 1.4 in the case $d=1$. We will use the following freeness 
result.

$\ $

{\bf 5.1 Proposition.} Let $\ncps$ be a $W^{*}$-probability space, with 
$\varphi$ faithful trace. Let $a,c$ be in $\A$, and assume that $c$ can
be factored as $c = up$, where $u \in \A$ is a unitary with 
Haar distribution, $p=p^{*} \in \A$ has a symmetric distribution, and
$\{ u,u^{*} \}$ is free from $\{ p \}$. (In other words, we assume that 
$c$ is $R$-diagonal.) If $\{ a,a^{*} \}$ is free from $\{ c,c^{*} \}$
in $\ncps$, then the selfadjoint matrices:
\begin{equation}
A \  = \left(  \begin{array}{cc}
0  &   a     \\
a^*  &  0
\end{array}   \right) , \ \ 
S \  = \left(  \begin{array}{cc}
0  &   c     \\
c^*  &  0
\end{array}   \right)  \ \ 
\end{equation}
are free in $\ncpst$.

$\ $

{\bf Proof.} We denote:
\[
{\cal X} \ := \ \{ u,u^{*} \} \ \cup \ 
\{ p^{k} - \varphi ( p^{k} ) I \ | \ k \geq 1 \} .
\]
A word made with letters from the alphabet $\X$ will be called 
``alternating'' if no two consecutive letters of the word are both
from $\{ u,u^{*} \}$ or both from $\{ p^{k} - \varphi ( p^{k} )I \ |$
$k \geq 1 \}$; the set of such alternating words will be denoted 
by $\Xalt$ Note that $\Xalt \subseteq$ Ker$( \varphi )$; this follows
(by using the definition of freeness) from the facts that $\X \subseteq$
Ker$( \varphi )$ and that $\{ u,u^{*} \}$ is free from $\{ p \}$.

Let us consider on the other hand the set:
\[
\Y \ = \ \Y_{11} \cup \Y_{12} \cup \Y_{21} \cup \Y_{22} , 
\]
where:
\[
\Y_{11} \ = \ \{ ( aa^{*} )^{k} - \varphi ( ( aa^{*} )^{k}) I) 
\ | \ k \geq 1 \} ,
\]
\[
\Y_{12} \ = \ \{ a ( a^{*}a )^{k} \ | \ k \geq 0 \} , \ \ 
\Y_{21} \ = \ \{ a^{*} ( aa^{*} )^{k} \ | \ k \geq 0 \} ,
\]
\[
\Y_{22} \ = \ \{ ( a^{*}a )^{k} - \varphi ( ( a^{*}a )^{k}) I)
\ | \ k \geq 1 \} .
\]
We will look at words of the form 
\begin{equation}
w \ = \ ( y_{1} - \lambda_{1} I ) x_{1}
( y_{2} - \lambda_{2} I ) x_{2} \cdots
( y_{n} - \lambda_{n} I ) x_{n} ,
\end{equation}
where $n \geq 1$, $y_{1}, \ldots , y_{n} \in \Y$, 
$\lambda_{1} , \ldots , \lambda_{n} \in \C$, $x_{1} , \ldots , x_{n}
\in \Xalt$, and where the following rules are obeyed:
\begin{equation}
\left\{  \begin{array}{l}
{ \mbox{if $y_{m} \in \Y_{11} \cup \Y_{21}$ ($1 \leq m \leq n$), then
$x_{m}$ begins with $u$;} }  \\
{ \mbox{if $y_{m} \in \Y_{12} \cup \Y_{22}$ ($1 \leq m \leq n$), then
$x_{m}$ begins with a $p^{k} - \varphi ( p^{k} ) I$;} }  \\
{ \mbox{if $y_{m} \in \Y_{11} \cup \Y_{12}$ ($2 \leq m \leq n$), then
$x_{m-1}$ ends with $u^{*}$;} }  \\
{ \mbox{if $y_{m} \in \Y_{21} \cup \Y_{22}$ ($2 \leq m \leq n$), then
$x_{m-1}$ ends with a $p^{k} - \varphi ( p^{k} ) I$;} }  \\
{ \mbox{if $y_{m} \in \Y_{11} \cup \Y_{22}$ ($1 \leq m \leq n$), then
$\lambda_{m} = 0$.} }
\end{array}  \right.
\end{equation}

We will prove the following:
\begin{equation}
\mbox{Claim: If $w$ satisfies (5.6), then $\varphi (w) = 0$.}
\end{equation}
The proof of the Claim (5.7) will be done by induction on the number 
$n$ of $x_{i}$'s and $y_{i}$'s entering the word $w$. For $n=1$, we have:
\[
\varphi (w) \ = \ \varphi ( \ (y_{1} - \lambda_{1}) x_{1} \ )
\]
\[
= \ \varphi ( y_{1} - \lambda_{1} ) \varphi ( x_{1} ) \ \
\mbox{ (because $\{ a,a^{*} \}$ free from $\{ c,c^{*} \}$) }
\]
\[
= \ 0 \ \ \mbox{ (because } \varphi ( x_{1} ) = 0).
\]

Let us next assume the Claim (5.7) is true for $n-1$, and prove it for
$n$. We first show that:
\begin{equation}
\varphi ( \ ( y_{1} - \lambda_{1} I ) x_{1}
( y_{2} - \lambda_{2} I ) x_{2} \cdots
( y_{n} - \lambda_{n} I ) x_{n} \ ) \ = 
\end{equation}
\[
\varphi ( \ ( y_{1} - \lambda_{1} ' I ) x_{1}
( y_{2} - \lambda_{2} ' I ) x_{2} \cdots
( y_{n} - \lambda_{n} ' I ) x_{n} \ ),
\]
for every $y_{1}, \ldots , y_{n} \in \Y$, 
$\lambda_{1} , \ldots , \lambda_{n},
\lambda_{1} ' , \ldots , \lambda_{n} ' \in \C$, 
$x_{1} , \ldots , x_{n} \in \Xalt$ such that the rules (5.6) are
satisfied. Clearly, it suffices to verify (5.8) in the situation when
$( \lambda_{1} , \ldots , \lambda_{n} )$ differs from 
$( \lambda_{1} ' , \ldots , \lambda_{n} ' )$ on only one position $k$, 
$1 \leq k \leq n$. For that $k$ we must have $y_{k} \in \Y_{12} \cup 
\Y_{21}$ (otherwise $\lambda_{k}$ and $\lambda_{k} '$ are both set to 
0 in (5.6)). But in such a case the difference of the two 
sides of (5.8) equals
\begin{equation}
( \lambda_{k} ' - \lambda_{k} )
\varphi ( \ ( y_{1} - \lambda_{1} I ) x_{1} \cdots 
( y_{k-1} - \lambda_{k-1} I ) x_{k-1} x_{k}
( y_{k+1} - \lambda_{k+1} I ) x_{k+1}  \cdots
( y_{n} - \lambda_{n} I ) x_{n} \ );
\end{equation}
and the quantity in (5.9) is indeed equal to 0, due to the induction
hypothesis. (The main point, in order to apply the induction hypothesis,
is to note that in both the possible cases -- $y_{k} \in \Y_{12}$,
$y_{k} \in \Y_{21}$ -- we will have $x_{k-1} x_{k} \in \Xalt$; this happens
because of the four ``concatenation'' rules stated in (5.6).)

Now, it is immediate that for every word $w = ( y_{1} - \lambda_{1} I ) 
x_{1} ( y_{2} - \lambda_{2} I ) x_{2} \cdots ( y_{n} - \lambda_{n} I ) 
x_{n}$ as in (5.5), we can find some new scalars $\lambda_{1} ' , \ldots ,
\lambda_{n} ' \in \C$ such that $w' = ( y_{1} - \lambda_{1} ' I ) 
x_{1} ( y_{2} - \lambda_{2} ' I ) x_{2} \cdots ( y_{n} - \lambda_{n} ' I ) 
x_{n}$ still satisfies the rules (5.6), and such that in addition:
\begin{equation}
\varphi ( y_{1} - \lambda_{1} ' I ) \ = \ \cdots \ = \ 
\varphi ( y_{n} - \lambda_{n} ' I ) \ = \ 0. 
\end{equation}
Indeed, if $1 \leq m \leq n$ is such that $y_{m} \in \Y_{12} \cup \Y_{21}$,
we can take $\lambda_{m} ' = \varphi ( y_{m} )$; while
if $1 \leq m \leq n$ is such that $y_{m} \in \Y_{11} \cup \Y_{22}$, then
the last rule (5.6) imposes $\lambda_{m} ' = 0 = \lambda_{m}$ -- but in
this case we also get $\varphi ( y_{m} ) = 0$ from the definitions of
$\Y_{11}, \Y_{22}$. The new word $w'$ satisfies $\varphi (w') = 0$; indeed,
besides (5.10) we also have $\varphi (x_{1}) = \cdots = \varphi (x_{n})
= 0$ (because $x_{1}, \ldots , x_{n} \in \Xalt$), and we only need to apply
the definition of freeness. Since (5.8) gives us that $\varphi (w) = 
\varphi (w')$, it follows that $\varphi (w) = 0$, and this concludes the 
proof of the Claim (5.7).

Let us finally look at the matrices $A,S$ defined in (5.4). In order to 
verify their freeness, it suffices to check that $\varphi_{2}  (W) = 0$
for every word:
\begin{equation}
W \ = \ \Bigl( A^{k_{1}} - \varphi_{2} ( A^{k_{1}} ) I_{2} \Bigr)
\Bigl( S^{l_{1}} - \varphi_{2} ( S^{l_{1}} ) I_{2} \Bigr) \cdots 
\Bigl( A^{k_{n}} - \varphi_{2} ( A^{k_{n}} ) I_{2} \Bigr)
\Bigl( S^{l_{n}} - \varphi_{2} ( S^{l_{n}} ) I_{2} \Bigr) ,
\end{equation}
with $n, k_{1}, l_{1}, \ldots , k_{n}, l_{n} \geq 1$. A straightforward
calculation shows that both diagonal entries of $W$ in (5.11) are words of
the type considered in (5.5)-(5.6). Hence the diagonal entries of $W$ are 
in Ker$( \varphi )$, by the Claim (5.7) -- and consequently $W \in$ 
Ker$( \varphi_{2} )$, as desired. {\bf QED}

$\ $

{\bf 5.2 Proposition} (the case $d=1$ of Theorem 1.4). Let $\nu$ be a 
probability measure with compact support on $[ 0, \infty )$, and let 
$\mu$ be the symmetric square root of $\nu$ (defined as in 3.1). Let 
$\ncps$ be a $W^{*}$-probability space, with $\varphi$ faithful trace,
and let $a \in \A$ be such that $a^{*}a$ has distribution $\nu$. Then:
\begin{equation}
\chi^{*} (a,a^{*}) \ \leq 2 \chi^{*} ( \mu ).
\end{equation}
Moreover, (5.12) holds with equality if $a$ is $R$-diagonal.

$\ $

{\bf Proof.} We may assume without loss of generality that there exists
a circular element $c \in \A$, of variance 1, such that $\{ c,c^{*} \}$
is free from $\{ a, a^{*} \}$. Then, by (5.1):
\begin{equation}
\chi^*(a,a^*)\ = \ \frac 12 \int_0^\infty\Bigl( \ \frac{2 }{1+t}
- \Phi^*(a+\sqrt t\,c, (a +\sqrt t\, c)^* )
\ \Bigr) dt+ \log(2\pi e).
\end{equation}

Consider on the other hand the space $\ncpst$ of $2 \times 2$-matrices
over $\ncps$, and the selfadjoint matrices $A,S \in M_{2} ( \A )$ defined
exactly as in Equation (5.4) of Proposition 5.1. Then $A$ has distribution
$\mu$ (by Remark 3.5), and is free from $S$ (by Proposition 5.1). 
From the form of $S$ it is immediate that
\[
\varphi_{2} ( S^{2n} ) = \varphi ( ( c^{*}c )^{n} ), \ \ 
\varphi_{2} (S^{2n+1}) = 0, \ \ \forall n \geq 0.
\]
It is known that $c^{*}c$ has the same distribution as the square of a 
semicircular element (see \cite{VDN}, Section 5.1); this implies that
$S$ is semicircular of variance 1.

Now, since $S$ is a normalized semicircular free from $A$, we can write:
\begin{equation}
\chi^* ( \mu ) \ = \  \chi^{*} (A) \ = \ \frac 12
\int_0^\infty\Bigl(\frac{1}{1+t} - \Phi^*(A+\sqrt{t}\, S)
\Bigr)dt + \frac{1}{2} \log(2\pi e).
\end{equation}
But for every $t \geq 0$:
\[
A + \sqrt{t} S \ = \ \left(  \begin{array}{cc}
0  &   a+ \sqrt{t} c   \\
(a + \sqrt{t} c)^*  &  0
\end{array}   \right) , \ \ 
\]
hence (again by Remark 3.5) the distribution of $A + \sqrt{t} S$ is the
symmetric square root of the distribution of
$(a + \sqrt{t} c)^{*} (a + \sqrt{t} c)$. When applied to this situation,
the Theorem 1.1 gives us that:
\begin{equation}
\Phi^{*} ( \ a + \sqrt{t} c, (a + \sqrt{t} c)^{*} \ ) \ \geq \
2 \Phi^{*} ( A + \sqrt{t} S ), \ \ \forall t \geq 0.
\end{equation}
The inequality (5.12) is obtained by replacing (5.15) in (5.13), and 
by comparing the result with (5.14).

If $a$ is $R$-diagonal, then so is $a+ \sqrt{t}c$ for every $t \geq 0$.
Indeed, $\sqrt{t}c$ is also $R$-diagonal, and the sum of two free 
$R$-diagonal elements is still $R$-diagonal (this follows for instance 
right away from the characterization of $R$-diagonality in terms of the
$R$-transform -- see \cite{NS}). But then the Theorem 1.1 implies that 
(5.15) holds with equality for every $t \geq 0$; and consequently, 
when we replace (5.15) in (5.13) and compare with (5.14), we obtain that 
(5.12) holds with equality too. {\bf QED}

$\ $

We now move to the proof of Theorem 1.5. We will use a known freeness
result, stated as follows.

$\ $

{\bf 5.3 Proposition.} Let $\ncps$ be a $W^{*}$-probability space,
with $\varphi$ faithful trace, let $\B \subseteq \A$ be a unital
$W^{*}$-subalgebra, and let $d$ be a positive integer. 

\vspace{4pt}

$1^{o}$ Let $\{ c_{ij} \}_{1 \leq i,j \leq d}$ be a family of
elements of $\A$ such that every $c_{ij}$ is circular of variance 
1, and such that $\{ c_{11}, c_{11}^{*} \} , \{ c_{12} ,
c_{12}^{*} \} , \ldots , \{ c_{dd}, c_{dd}^{*} \} , \B$
are free. Then the matrix $C = \matcij$ is a circular element of 
variance $d$ in $M_{d} ( \A )$, and $\{ C,C^{*} \}$ is free 
from $M_{d} ( \B )$.

\vspace{4pt}

$2^{o}$ Let $\{ s_{ij} \}_{1 \leq i,j \leq d}$ be a family of
elements of $\A$ such that: $s_{ij}^{*} = s_{ji}$, for every 
$1 \leq i,j \leq d$; $s_{ii}$ is semicircular of variance 1 for 
every $1 \leq i \leq d$; $s_{ij}$ is circular of variance 1 for 
every $1 \leq i<j \leq d$; and $\{ s_{11} \} , \ldots , 
\{ s_{dd} \}$, $\{ s_{12}, s_{12}^{*} \} , \ldots , \{ s_{d-1,d} 
s_{d-1,d}^{*} \} , \B$ are free. Then the selfadjoint matrix
$S =  ( s_{ij} )_{i,j=1}^{d}$ is a semicircular element of 
variance $d$ in $M_{d} ( \A )$, and is free from $M_{d} ( \B )$.

$\ $

For the fact that $C$ of $5.3.1^{o}$ is circular and that $S$ of
$5.3.2^{o}$ is semicircular, see \cite{VDN}, Section 5.1; for the 
additional assertions concerning the freeness from $M_{d} ( \B )$,
see \cite{Sh}. For the sake of completeness, we indicate a way of 
proving Proposition 5.3 which, quite amusingly, comes out directly 
from the considerations of the preceding section. We have:

$\ $

{\bf 5.4 Lemma.} Let $\ncps$ be a $W^{*}$-probability space, with 
$\varphi$ faithful trace. Let $\{ s_{i} \}_{1 \leq i \leq k}$ be a 
family of selfadjoint elements of $\A$, and let $\B \subseteq \A$
be a unital $W^{*}$-subalgebra. Assume that 
$\{ s_{i} \}_{1 \leq i \leq k}$ is its own conjugate with respect 
to $\B$. Then every $s_{i}$ ($1 \leq i \leq k$) is semicircular of
variance 1, and $\{ s_{1} \}, \ldots , \{ s_{k} \} , \B$ are free.

$\ $

{\bf Proof of Lemma 5.4.} This is an immediate consequence of the 
free Cramer-Rao inequality, as stated in \cite{V-P5}, Proposition 6.9.
The line of the argument goes as follows. By enlarging $\ncps$ if 
necessary, we can assume that there also exists in $\A$ a family 
$\{ s_{i} ' \}_{1 \leq i \leq k}$ of semicircular elements of variance 
1, such that $\{ s_{1} ' \}, \ldots , \{ s_{k} ' \}$, $\B$ are free. 
Then $\{ s_{i} ' \}_{1 \leq i \leq k}$ is its own conjugate with respect 
to $\B$, by Propositions 3.8 and 3.6 of \cite{V-P5}. The hypothesis
that $\{ s_{i} \}_{1 \leq i \leq k}$ is its own conjugate with respect 
to $\B$ amounts to the fact that
\[
\varphi ( s_{i} b_{0} s_{i_{1}} b_{1} \cdots s_{i_{n}} b_{n} ) \ = \ 
\sum_{m=1}^{n} \delta_{i,i_{m}} 
\varphi ( b_{0} s_{i_{1}} \cdots s_{i_{m-1}} b_{m-1} ) 
\cdot \varphi ( b_{m} s_{i_{m+1}} \cdots s_{i_{n}} b_{n} )
\]
\begin{equation}
= \ \sum_{m=1}^{n} \delta_{i,i_{m}} 
\varphi \Bigl( s_{i_{1}} b_{1} \cdots s_{i_{m-1}} ( b_{m-1} b_{0} ) \Bigr) 
\cdot \varphi \Bigl( s_{i_{m+1}} b_{m+1} \cdots s_{i_{n}}
( b_{n} b_{m} ) \Bigr) ,
\end{equation}
for every $n \geq 0$, $b_{0}, \ldots , b_{n} \in \B$,
$1 \leq i,i_{1}, \ldots , i_{n} \leq k$. Since 
$\{ s_{i} ' \}_{1 \leq i \leq k}$ also has the property of being its own 
conjugate with respect to $\B$, (5.16) remains true when we replace $s_{i}$
by $s_{i} '$ and $s_{i_{1}}$ by $s_{i_{1}} ' , \ldots , s_{i_{n}}$ by
$s_{i_{n}} '$. But then an induction argument immediately gives that: 
\begin{equation}
\varphi ( s_{i} b_{0} s_{i_{1}} b_{1} \cdots s_{i_{n}} b_{n} ) \ = \ 
\varphi ( s_{i} ' b_{0} s_{i_{1}} ' b_{1} \cdots s_{i_{n}} ' b_{n} ),
\end{equation}
for every $n \geq 0, \ b_{0}, \ldots , b_{n} \in \B , \
1 \leq i,i_{1}, \ldots , i_{n} \leq k$.
Finally, from (5.17) and the fact that $s_{1} ', \ldots s_{k} '$ are 
normalized semicirculars, with 
$\{ s_{1} ' \}, \ldots , \{ s_{k} ' \} , \B$ free, it follows that 
$s_{1} , \ldots , s_{k}$ also have these properties. {\bf QED}

$\ $

{\bf Proof of Proposition 5.3.} The proofs of $1^{o}$ and $2^{o}$
are similar; we will show $1^{o}$, and leave $2^{o}$ as an exercise
to the reader.

By working with the real and imaginary parts of the elements $c_{ij}$,
and by using Propositions 3.8 and 3.6 of \cite{V-P5}, one obtains that 
the conjugate of $\cijcijs$ with respect to $\B$ is $\xiijxiijs$,
with $\xi_{ij} := c_{ij}^{*}$, $1 \leq i,j \leq d$. The Proposition 4.1
from the preceding section (or rather its proof) applies to this 
situation, and gives that $\{ X,X^{*} \}$ is the conjugate of 
$\{ C,C^{*} \}$ with respect to $M_{d} ( \B )$, where:
\begin{equation}
X \ := \ \frac{1}{d} ( \xi_{ji} )_{i,j =1}^{d} \ = \ 
\frac{1}{d} C^{*} .
\end{equation}
From (5.18) it is immediate that if we set $S_{1} = (C+C^{*})/ \sqrt{2d}$,
$S_{2} = (C-C^{*})/i \sqrt{2d}$, then $\{ S_{1} , S_{2} \}$ is its own 
conjugate with respect to $M_{d} ( \B )$. But then we can use the Lemma 5.4
to infer that $S_{1},S_{2}$ are semicirculars of variance 1 in 
$M_{d} ( \A )$, such that $\{ S_{1} \} , \{ S_{2} \} , M_{d} ( \B )$ are 
free. This in turn implies that $C = \sqrt{d/2} ( S_{1} + i S_{2} )$ is
circular of variance $d$, and free from $M_{d} ( \B )$. {\bf QED}

$\ $

{\bf Proof of Theorem 1.5.} The proofs of $1.5.1^{o}$ and $1.5.2^{o}$
are similar; in order to offer the reader a variation, we will this
time show $2^{o}$, and leave $1^{o}$ as an exercise.

The selfadjoint family $\bijs$ appearing on the left-hand side of 
(1.10) is to be looked at as
$\{ b_{ij}, b_{ij}^{*} \}_{1 \leq i<j \leq d} \cup 
\{ b_{ii} \}_{1 \leq i \leq d}$; thus (5.1) applies with 
$m = d(d-1)/2$, $n=d$, and yields the formula:
\begin{equation}
\chi^*(\ \bijs\ )\ = \ \frac 12 \int_0^\infty\Bigl(\frac{d^2
}{1+t}-\Phi^*(\ \{b_{ij}+\sqrt{t}\, s_{ij}
\}_{1\leq i,j\leq d}\ )\Bigr)dt+\frac{d^2}2 \log(2\pi e),
\end{equation}
where the family $\{ s_{ij} \}_{1 \leq i,j \leq d}$ of elements
of $\A$ has the following 
properties: $s_{ij}^{*} = s_{ji}$, for every $1 \leq i,j \leq d$;
$s_{ii}$ is semicircular of variance 1 for every $1 \leq i \leq d$;
$s_{ij}$ is circular of variance 1 for every $1 \leq i<j \leq d$;
and the sets $\{ s_{11} \} , \ldots , \{ s_{dd} \}$,
$\{ s_{12}, s_{12}^{*} \} , \ldots , \{ s_{d-1,d} s_{d-1,d}^{*} \}$,
$\{ b_{ij} \ | \ 1 \leq i,j \leq d \}$ are free.

If we denote $S := ( s_{ij} )_{i,j=1}^{d} \in M_{d} ( \A )$, then
$d^{-1/2}S$ is semicircular of variance 1, free from $B$ in $\ncpsd$
(by Proposition $5.3.2^{o}$, where the choice of $W^{*}$-subalgebra
$\B \subseteq \A$ is made to be $\B := W^{*} ( \{ I \} \cup  \bijs )$ ). 
We can therefore use $d^{-1/2}S$ in the calculation of the free entropy 
$\chi^{*} (B)$; it is in fact more convenient to write the formula
for $\chi^{*} ( d^{-1/2}B )$:
\begin{equation}
\chi^{*} ( \frac{1}{\sqrt{d}} B ) \ = \ \frac 12
\int_0^\infty\Bigl(\frac{1}{1+t} - \Phi^*( \frac{1}{\sqrt{d}} B  + 
\sqrt{\frac{t}{d}} S) \Bigr)dt + \frac{1}2 \log(2\pi e) .
\end{equation}
The scaling formulas for $\Phi^{*}$ and $\chi^{*}$ are
\[
\Phi^{*} ( \lambda x ) = \lambda^{-2} \Phi^*(x), \ \ 
\chi^{*} ( \lambda x ) = \chi^*(x) + \log ( \lambda )
\]
(for $\lambda > 0$ and $x$ a selfadjoint random variable -- see 
\cite{V-P5}, Sections 6.2(b) and 7.8). Thus (5.20) can also be written 
in the form:
\begin{equation}
\chi^{*} ( B ) - \frac{\log d}{2}  \ = \ \frac 12
\int_0^\infty\Bigl(\frac{1}{1+t} - d \Phi^*( B  + \sqrt{t} S) \Bigr)dt
+ \frac{1}2 \log(2\pi e) .
\end{equation}

Now, for every $t \geq 0$, the Theorem $1.2.2^{o}$ gives us:
\begin{equation}
\Phi^{*} ( \ \{ b_{ij} + \sqrt{t} s_{ij} \}_{1 \leq i,j \leq d} \ ) 
\ \geq \ d^{3} \Phi^*( B  + \sqrt{t} S).
\end{equation}
If we replace (5.22) into (5.19), and compare the result with (5.21), 
then (1.10) of Theorem 1.5 is obtained.

If $B$ is free from the algebra of scalar matrices $M_{d} ( \C I )
\subseteq M_{d} ( \A )$, then the same is true for $B + \sqrt{t} S$, 
for every $t \geq 0$; this is because (as implied by Proposition 
$5.3.2^{o}$) $S$ is free from $M_{d} ( W^{*} ( \bijs ) )$, which in 
turn implies that $\{ B,S \}$ is free from $M_{d} ( \C I )$. But in 
this situation, the Theorem $1.2.2^{o}$ implies that (5.22) holds with
equality, for every $t \geq 0$; and the same argument used in the 
preceding paragraph shows now that (1.10) holds with equality, too.
{\bf QED}

$\ $

{\bf Proof of Theorem 1.4.} This follows from Proposition 5.2 and 
Theorem $1.5.1^{o}$, exactly as in the same way as Theorem 1.3 was
obtained from 1.1 and $1.2.1^{o}$ at the end of Section 4.
{\bf QED}

$\ $
 
$\ $
 
$\ $

\end{document}